\documentclass[10pt,a4paper,reqno]{amsart}
\usepackage{amsfonts,amsthm,latexsym,amsmath,amssymb,amscd,amsmath}
\usepackage{epsfig,psfrag}

\newtheorem{theo}{Theorem}
\newtheorem{fact}{Fact}
\newtheorem{prob}{Problem}
\newtheorem{prop}{Proposition}[section]
\newtheorem{coro}[prop]{Corollary}

\newtheorem{lemma}[prop]{Lemma}

\newcommand{\chop}{{\mathbf{chop}}}
\newcommand{\Diag}{{\mathrm{Diag}}}
\newcommand{\AD}{{\mathbf{AD}}}
\newcommand{\TR}{{\mathbf{TR}}}
\newcommand{\NE}{{\mathbf{NE}}}
\newcommand{\NW}{{\mathbf{NW}}}
\newcommand{\SE}{{\mathbf{SE}}}
\newcommand{\SW}{{\mathbf{SW}}}
\newcommand{\Bi}{{\mathbf{B}}}
\newcommand{\Ci}{{\mathbf{C}}}
\newcommand{\bU}{{\mathbf{U}}}
\newcommand{\aux}{{\mathbf{aux}}}

\newcommand{\ba}{{\bf a}}

\newcommand{\Bruhat}{\operatorname{Bru}}

\newcommand{\trace}{\operatorname{trace}}
\newcommand{\trd}{\operatorname{trd}}
\newcommand{\diag}{\operatorname{diag}}
\newcommand{\Spin}{\operatorname{Spin}}

\newcommand{\ZZ}{{\mathbb{Z}}}

\newcommand{\RR}{{\mathbb{R}}}

\newcommand{\Ss}{{\mathbb{S}}}
\newcommand{\BB}{{\mathbb{B}}}
\newcommand{\PP}{{\mathbb{P}}}

\newcommand{\ga}{\gamma}

\newcommand{\F}{\mathfrak F}
\newcommand{\tF}{\tilde\F}

\newcommand{\fT}{\mathfrak T}
\newcommand{\1}{{\mathbf 1}}
\newcommand{\U}{{\mathcal U}}
\newcommand{\Little}{{\mathcal L}}
\newcommand{\fast}{{\star}}

\newcommand{\nobf}{\noindent\bf}

\def\qed{\unskip\nobreak\hfil\penalty50\hskip1.75em\null\nobreak\hfil
$\blacksquare$ {\parfillskip=0pt \finalhyphendemerits=0 \par}\goodbreak}

\begin{document}
\title[Locally convex curves in $\Ss^n$ and  the group $B^{+}_{n+1}$]
{Spaces of locally convex curves in $\Ss^n$
and combinatorics of the group $B^{+}_{n+1}$}

\author{Nicolau C. Saldanha and Boris Shapiro}
\dedicatory {Dedicated  to the memory of the one and only
Vladimir Igorevich Arnold}

\begin{abstract}
In the 1920's Marston Morse developed what is now known as Morse theory
trying to study the topology of the space of closed curves
on $\Ss^2$ (\cite{Morse}, \cite {Klingenberg}).
We propose to attack a very similar problem,
which 80 years later remains open,
about the topology of the space of closed curves on $\Ss^2$ 
which are locally convex (i.e., without inflection points).
One of the main difficulties is
the absence of the covering homotopy principle for the map sending a
non-closed locally convex curve to the Frenet frame at its endpoint.

In the present paper we study the spaces of locally
convex curves  in $\Ss^n$ with a given initial and final Frenet frames.
Using combinatorics of $B^{+}_{n+1} = B_{n+1} \cap SO_{n+1}$,
where $B_{n+1} \subset O_{n+1}$ is the usual Coxeter-Weyl group,
we show that for any $n\ge 2$ these spaces fall in at most
$\lceil\frac{n}{2}\rceil+1$ equivalence classes up to homeomorphism.
We also study this classification in the double cover $\Spin(n+1)$.
For $n = 2$ our results complete the classification
of the corresponding spaces into two topologically distinct classes,
or three classes in the spin case.
\end{abstract}

\maketitle


\medbreak

\section{Introduction and main results} 


In what follows we will study different spaces of curves
$\gamma:[0,1]\to \Ss^n$ (or to $\RR^{n+1}$);
we start with some basic definitions. 
A smooth curve $\gamma:[0,1]\to\RR^{n+1}$ is called \textit{locally convex}
if its Wronskian
$W_\gamma(t)=\det \left(\ga(t), \ga'(t),\ga''(t),\dots, \ga^{(n)}(t)\right)$
is non-vanishing for all $t\in[0,1]$
(see \cite{Anisov}, \cite{SK}, \cite{Shapiro2}, \cite{ShapiroM}). 
A smooth curve $\gamma$ is called \textit{(globally) convex}
if for any linear hyperplane $H\subset \RR^{n+1}$ 
the intersection $H\cap \ga$ consists of at most $n$
points counting  multiplicities;  
it is an easy exercise to  check that global convexity implies local.
Observe that if $\gamma: [0,1] \to \RR^{n+1}$ is locally convex
then so is its spherical projection
$\ga/|\ga|: [0,1]\to \Ss^n \subset \RR^{n+1}$.
Notice that for $n = 2$, a curve $\gamma:[0,1]\to \Ss^2$
is locally convex if its geodesic curvature is never zero
(and therefore has constant sign) and a closed curve $\gamma:[0,1]\to \Ss^2$
is globally convex if it is the boundary of the intersection
of the sphere with a convex cone.

For various technical reasons,
the space of smooth curves is too small and not the most adequate.
The definition of local convexity makes sense for other spaces,
such as the Banach spaces $C^r$, $r \ge n$
and the Sobolev spaces $H^r$, $r > n$.
In Section 2 below we shall introduce an ``official'' topology
for the spaces of locally convex curves:
this turns out to be a Hilbert space containing
all the above spaces.
As with other questions concerning infinite dimensional topology,
the choice of space actually has little consequence.

Locally convex curves in $\RR^{n+1}$ are closely related
to fundamental solutions of   linear 
ordinary  homogeneous differential equations of order $n+1$
on $[0,1]$ with real-valued coefficients.
Namely,  
if $y_0, y_1, \ldots, y_n$ are linearly independent solutions
of an equation 
\[ y^{(n+1)} + a_n(t) y^{(n)} + \cdots + a_0(t) y = 0 \]
with $a_i(t)\in C[0,1]$, $i=0,\ldots,n$ ,
then $\gamma = (y_0, y_1, \ldots, y_n)$ is locally convex.
A locally convex $\ga$ is called \textit {positive} if $W_\ga(t)>0$
and \textit{negative} otherwise.
From now on we mostly consider positive curves.  

Given a smooth positive locally convex $\gamma:[0,1]\to\RR^{n+1}$,
define its Frenet frame $\F_\ga: [0,1] \to SO_{n+1}$
as the result of the Gram-Schmidt orthogonalization
of its Wronski curve
$\left(\ga(t), \ga'(t),\ga''(t), \dots, \ga^{(n)}(t)\right)$.
In other words, $\F_\ga$ satisfies the relation 
$\left(\ga(t), \ga'(t),\ga''(t), \dots, \ga^{(n)}(t)\right)
= \F_\ga(t) R(t)$ where $R(t)$ is an upper triangular matrix
with positive diagonal.
Let $\Little\Ss^n$ be  the space of all positive locally convex curves
$\gamma:[0,1] \to \Ss^n$
(in the appropriate space, to be defined in Section 2)
with the standard initial frame $\F_\ga(0) = I$,
where $I \in SO_{n+1}$ is the identity matrix of size $(n+1)\times (n+1)$.
As we shall see (Lemma \ref{lemma:Lcontract}),
the space $\Little\Ss^n$ is a contractible Hilbert manifold
and therefore diffeomorphic to Hilbert space.




Given $Q \in SO_{n+1}$, let $\Little\Ss^n(Q) \subset \Little\Ss^n$
be the set of positive locally convex curves on $\Ss^n$ with the 
standard initial and the prescribed final frame $\F_\ga(1) = Q$;
one of the main difficulties is that the map $\Little\Ss^n \to SO_{n+1}$
taking $\gamma$ to $\F_\ga(1)$ is not a fibre bundle.
Let $\Pi: \Spin_{n+1} \to SO_{n+1}$ ($n\ge 2$)
be the universal cover (recall that this cover is a double cover).
Denote by $\1 \in \Spin_{n+1}$ the identity element 
and by  $-\1 \in \Spin_{n+1}$
the unique nontrivial element with $\Pi(-\1) = I$.
For $\gamma \in \Little\Ss^n$, the map $\F_\ga:[0,1]\to SO_{n+1}$
can be uniquely lifted to $\tF_\ga:[0,1] \to \Spin_{n+1}$,
$\F_\gamma = \Pi \circ \tF_\gamma$, $\tF_\gamma(0) = \1$.
Given $z \in \Spin_{n+1}$, let $\Little\Ss^n(z) \subset \Little\Ss^n(\Pi(z))$
be the set of positive locally convex curves $\gamma \in \Little\Ss^n(\Pi(z))$
with $\tF_\gamma(1) = z$. 
One can immediately observe that
$\Little\Ss^n(\Pi(z)) = \Little\Ss^n(z) \sqcup \Little\Ss^n(-z)$.
The Hilbert manifolds $\Little\Ss^n(Q)$ and $\Little\Ss^n(z)$ for various 
$Q \in SO_{n+1}$ and $z \in \Spin_{n+1}$ 
are the main objects of study in this paper.

Some information  about the topology of $\Little\Ss^n(Q)$,
mostly in the case $Q=I$ or  in the case $n=2$,
was earlier obtained in \cite{Anisov}, \cite{Little},
\cite{Saldanha1}, \cite{Saldanha2}, \cite{SK},
\cite{Shapiro2} and \cite{ShapiroM}.
In particular, it was shown that the number of connected components of
$\Little\Ss^n(I)$ equals $3$ for even $n$
and $2$ for odd $n>1$, which  is related to the existence of closed globally convex 
curves on all even-dimensional spheres.
It was also shown in \cite{Anisov}
that for $n$ even the space of closed globally convex curves with a fixed
initial frame is contractible.
The first nontrivial information about the
higher homology and homotopy groups of these components
can be found in \cite{Saldanha1} and \cite{Saldanha2}. 

In this paper we leave aside the fascinating and widely open question 
about  the topology of the spaces $\Little\Ss^n(Q)$ and
concentrate on the following. 

\begin{prob}
\label{prob:main}
How many different (i.e., non-homeomorphic) spaces 
are there among $\Little\Ss^n(Q)$,  $Q \in SO_{n+1}$, $n\ge 2$?
Analogously, how many different spaces are there among
$\Little\Ss^n(z)$,  $z \in Spin_{n+1}$?
\end{prob}

To formulate our partial answer to the latter question
we need to introduce the following set of matrices. 
For a positive integer  $m$ let
\[ M^m_s = \diag(-1,\ldots,-1,1,\ldots,1), \qquad
s \in \ZZ, \quad |s| \le m, \quad s \equiv m \pmod 2,
 \]
be the diagonal $m \times m$ matrix whose 
first $(m-s)/2$ entries equal to $-1$ and the remaining
$(m+s)/2$ entries equal to $1$.
Notice that $s$ equals both the trace and the signature of $M^m_s$
and that $M^m_s \in SO_m$ if and only if $s \equiv m \pmod 4$. 
In the latter case, let $\pm w^m_s \in \Spin_m$
be the two preimages of $M^m_s \in SO_m$.


Our first result is as follows.

\begin{theo}
\label{theo:main}
For  $n \ge 2$, any $Q \in SO_{n+1}$,
and any $z \in \Spin_{n+1}$  one has:  
\begin{enumerate}
\item{Each space $\Little\Ss^n(Q)$ is homeomorphic  to one of the subspaces
$\Little\Ss^n({M^{n+1}_s})$, where $|s| \le n+1$, $s \equiv n+1 \pmod 4$
(there are $\lceil\frac{n}{2}\rceil+1$ such subspaces).}
\item{For $n$ even, each space $\Little\Ss^n(z)$, $z \in \Spin_{n+1}$,
is homeomorphic to one of the subspaces $\Little\Ss^n(\1)$,
$\Little\Ss^n(-\1)$,
$\Little\Ss^n(w^{n+1}_{n-3})$,
$\Little\Ss^n(w^{n+1}_{n-7})$, \dots, 
$\Little\Ss^n(w^{n+1}_{-n+5})$,
$\Little\Ss^n(w^{n+1}_{-n+1})$.  }
\item{For $n$ odd, each space $\Little\Ss^n(z)$, $z \in \Spin_{n+1}$,
is homeomorphic  to one of the subspaces
$\Little\Ss^n(\1)$, $\Little\Ss^n(-\1)$,
$\Little\Ss^n(w^{n+1}_{n-3})$,
$\Little\Ss^n(w^{n+1}_{n-7})$, \dots, 
$\Little\Ss^n(w^{n+1}_{-n+3})$,
$\Little\Ss^n(w^{n+1}_{-n-1})$, $\Little\Ss^n(-w^{n+1}_{-n-1})$. }
\end{enumerate} 
\end{theo} 

Using Theorem \ref{theo:dois} below and results proved elsewhere
(\cite{Saldanha1}, \cite{Saldanha2}) we check that for $n=2$ 
the above spaces are pairwise non-homeomorphic. 
It is natural to ask whether
they are likewise non-homeomorphic for $n \ge 3$;
see discussions in the first subsection of the conclusion. 

We might want to describe the topology of these spaces;
the next result gives some partial answers.
Let $\Omega SO_{n+1}(Q)$ (resp. $\Omega\Spin_{n+1}(z)$)
be the space of all continuous curves
$\alpha: [0,1] \to SO_{n+1}$ (resp. $\alpha: [0,1] \to \Spin_{n+1}$)
with $\alpha(0) =  I$ and $\alpha(1) = Q$
(resp. $\alpha(0)=\1$ and $\alpha(1) =z$).
Using the Frenet frame we define \textit{Frenet frame injections}:
\[
\begin{array}{rcl}
\F_{[Q]}: \Little\Ss^n(Q) & \to & \Omega SO_{n+1}(Q), \\
\gamma &\mapsto & \F_\gamma
\end{array}
\qquad
\begin{array}{rcl}
\tilde\F_{[z]}: \Little\Ss^n(z) & \to & \Omega\Spin_{n+1}(z). \\
\gamma & \mapsto & \tilde\F_\gamma 
\end{array}
\]
It is a classical fact that the value of $Q$ (resp.  $z$) 
does not change the space $\Omega SO_{n+1}(Q)$
(resp. $\Omega\Spin_{n+1}(z)$) up to homeomorphism.  
Therefore, we usually omit $Q$ (resp. $z$) and  write $\Omega SO_{n+1}$
(resp. $\Omega\Spin_{n+1}$) instead. 


\begin{theo}
\label{theo:dois}
For $n \ge 2$, consider the Frenet frame injections as above.
\begin{enumerate}
\item{For all $Q \in SO_{n+1}$ and for all $z \in \Spin_{n+1}$
the maps $\F_{[Q]}$ and $\tilde\F_{[z]}$
are weakly homotopically surjective.}
\item{If $|s| \le 1$ then the Frenet frame injections
$\F_{[{M^{n+1}_s}]}$ and $\tilde\F_{[{w^{n+1}_s}]}$
are weak homotopy equivalences.
In this case there exist homeomorphisms
\[ \Little\Ss^n({M^{n+1}_s}) \approx \Omega SO_{n+1}, \qquad
\Little\Ss^n({w^{n+1}_s}) \approx \Omega\Spin_{n+1}. \]}
\end{enumerate}
\end{theo} 

Recall that a map $X \to Y$ is  weakly homotopically surjective
if the induced maps $\pi_k(X) \to \pi_k(Y)$ are surjective;
also, a map $X \to Y$ is a weak homotopy equivalence
if the induced maps $\pi_k(X) \to \pi_k(Y)$ are isomorphisms.

Notice that, in general,  
for arbitrary $Q$ or $z$  it is by no means true
that the Frenet frame injection induces a homotopy equivalence:
even the number of connected components can be different.


Versions of Theorems \ref{theo:main} and \ref{theo:dois}
also hold for the spaces $C^k \cap \Little\Ss^n(Q)$ and
$C^k \cap \Little\Ss^n(z)$.
These facts follow from our results
together with Theorem 2 in \cite{BST};
alternatively, our proofs can be adapted
(with some extra rather routine work).

\noindent 
\textit {Acknowledgements.}
The first named author gratefully acknowledges the support
of CNPq, FAPERJ and CAPES (Brazil) and is sincerely grateful to
the mathematics department of the Stockholm University for
the kind hospitality during his visits to Sweden in 2005 and 2007.  

\bigbreak

\section{Frenet frames and Jacobian curves}

We collect in this section a few basic notions and  facts.
The logarithmic derivative of a curve $\Gamma: [0,1] \to SO_{n+1}$
is defined as \( \Lambda(t) = (\Gamma(t))^{-1} \Gamma'(t) \).
Notice that \( \Lambda(t) \)
belongs to the Lie algebra and is therefore
automatically skew symmetric.
Let $\fT \subset so_{n+1}$ be the set of
tridiagonal skew symmetric matrices
with positive subdiagonal entries,
or skew Jacobi matrices, i.e.,
of matrices of the form
\[ \begin{pmatrix}
    & -c_1 &        &        &  \\
c_1 &      & -c_2   &        &  \\
    &  c_2 &        & \ddots &  \\
    &      & \ddots &        & -c_n \\
    &      &        &  c_n   & 
\end{pmatrix}, \qquad
c_i > 0. \]
A curve  $\Gamma: I \to SO_{n+1}$ is called \textit{Jacobian}
if its logarithmic derivative $\Lambda$ satisfies
$\Lambda(t) \in \fT$ for all $t \in I$
(where $I \subset \RR$ is an interval).

\begin{lemma}
\label{lemma:Frenet}
Let $\Gamma: [0,1] \to SO_{n+1}$ be a smooth curve with $\Gamma(0) = I$.
The curve $\Gamma$ is Jacobian if and only if
there exists $\gamma \in \Little\Ss^n$ with $\F_\gamma = \Gamma$.
\end{lemma}

Recall that a smooth curve $\gamma: [0,1] \to \Ss^n$
belongs to $\Little\Ss^n$ if and only if $\F_\gamma(0) = I$
and $\gamma$ is (positive) locally convex:
\[ \det \left(\gamma(t), \gamma'(t),\gamma''(t),\ldots,
\gamma^{(n)}(t) \right) > 0. \]

\begin{proof}
Consider $\gamma \in \Little\Ss^n$ and its Wronski curve
\[ G(t) = \left( \gamma(t) \; \gamma'(t) \; \cdots
\; \gamma^{(n)}(t) \right)
=  \F_\ga(t) R(t). \]
We have 
\[ G'(t) = \left( \gamma'(t) \; \gamma''(t) \; \cdots
\; \gamma^{(n+1)}(t) \right)
= G(t) H(t) \]
for $H(t)$ an upper Hessenberg matrix whose subdiagonal entries equal to $1$:
\[ H(t) = 
\begin{pmatrix}
0 & 0 & 0 & \cdots & 0 & \ast \\
1 & 0 & 0 & \cdots & 0 & \ast \\
0 & 1 & 0 & \cdots & 0 & \ast \\
& \vdots & & \vdots & \\
0 & 0 & 0 & \cdots & 1 & \ast
\end{pmatrix}. \]
Recall that $H$ is upper Hessenberg if $(H)_{ij} = 0$ whenever $i > j+1$.
Write $\Gamma = \F_\ga$ and substitute $\Gamma(t) R(t)$ for $G(t)$
in the equations above to obtain
\[ \Gamma'(t) R(t) + \Gamma(t) R'(t) = \Gamma(t) R(t) H(t) \]
and therefore
\[  \Lambda(t) = (\Gamma(t))^{-1} \Gamma'(t) =
- R'(t) (R(t))^{-1} + R(t) H(t)(R(t))^{-1} \]
which is upper Hessenberg with positive subdiagonal entries
(the first product is upper triangular, 
the second one is upper Hessenberg).
Since we know that $\Lambda(t) \in so_{n+1}$,
we have $\Lambda(t) \in \fT$, proving one implication.

For the other implication, consider $\Gamma: [0,1] \to SO_{n+1}$
such that $\Gamma(0) = I$, $\Gamma'(t) = \Gamma(t) \Lambda(t)$
and \( \Lambda(t) \in \fT \) for all $t \in [0,1]$.
Set $\gamma(t) = \Gamma(t) e_1$.
We have $\gamma'(t) = \Gamma'(t) e_1 = \Gamma(t) \Lambda(t) e_1 
= \Gamma(t) (\Lambda(t))_{21} e_2 = p_1(t) \Gamma(t) e_2$,
$p_1(t) > 0$.
Similarly,
\[ \gamma''(t) = (p_1)'(t) \Gamma(t) e_2 + p_1(t) \Gamma(t) \Lambda(t) e_2 =
p_2(t) \Gamma(t) e_3 + r_{22}(t) \Gamma(t) e_2 + r_{21}(t) \Gamma(t) e_1, \]
where $p_2(t) > 0$ and the values of $r_{ij}(t)$ are not important.
In general
\[ \gamma^{(j)}(t) = p_j(t) \Gamma(t) e_{j+1} +
\sum_{i \le j} r_{j i}(t) \Gamma(t) e_{i}, \qquad
p_j(t) > 0. \]
Thus applying Gram-Schmidt to the Wronski curve
\[ \left( \gamma(t)\; \gamma'(t)\; \cdots\; \gamma^{(n)}(t) \right) \]
yields $\F_\gamma(t) = \Gamma(t)$,
completing the proof.
\end{proof}

A smooth Jacobian curve $\Gamma: I \to SO_{n+1}$ is called 
\textit{globally Jacobian} if $\gamma: I \to \Ss^n$,
$\gamma(t) = \Gamma(t) e_1$, is globally convex.

Notice that given a smooth function $\Lambda: [0,1] \to \fT \subset so_{n+1}$
the initial value problem
\[ \Gamma'(t) = \Gamma(t) \Lambda(t), \quad \Gamma(0) = I \eqno{(\ast)} \]
yields $\Gamma$ as in the lemma and therefore
a smooth curve $\gamma \in \Little\Ss^n$.
This establishes a homeomorphism between the space
of smooth curves $\gamma \in \Little\Ss^n$
and the convex set of smooth functions $\Lambda: [0,1] \to \fT$.
We will denote this correspondence by
\[ \Lambda_\gamma(t) = (\F_\gamma(t))^{-1} (\F_\gamma)'(t). \]
It will be convenient to have examples of locally convex curves
and corresponding Jacobian curves.

\begin{lemma}
\label{lemma:phonewire}
For  $n + 1 = 2k$ let
$c_i$ and $a_i$ ($i = 1, \ldots, k$) be positive parameters with
$a_i$ mutually distinct and $c_1^2 + \cdots + c_k^2 = 1$. Set
\[ \xi(t) = \left( c_1 \cos(a_1 t), c_1 \sin(a_1 t), \ldots,
c_k \cos(a_k t), c_k \sin(a_k t) \right). \]
For  $n + 1 = 2k + 1$ let
$c_0$, $c_i$ and $a_i$ ($i = 1, \ldots, k$) be positive parameters with
$a_i$ mutually distinct and $c_0^2 + c_1^2 + \cdots + c_k^2 = 1$. Set
\[ \xi(t) = \left( c_0, c_1 \cos(a_1 t), c_1 \sin(a_1 t), \ldots,
c_k \cos(a_k t), c_k \sin(a_k t) \right). \]
In both cases the curve $\xi: [0,1] \to \Ss^n$ is locally convex
with constant $\Lambda_\xi$.
Conversely, if $\tilde\xi: [0,1] \to \Ss^n$ is locally convex 
with $\Lambda_{\tilde\xi}$ constant then $\tilde\xi = Q\xi$
for some $Q \in SO_{n+1}$ and $\xi$ as above
(for appropriate $c_i$ and $a_i$).  
Furthermore,
assume $a_i/(4\pi) \in \ZZ$ and set $Q = (\F_\xi(0))^{-1}$
and $\xi_1(t) = Q\xi(t)$: we have $\xi_1 \in \Little\Ss^n(\1)$.
\end{lemma}

\begin{figure}[ht]
\begin{center}
\psfrag{(a)}{(a)}
\psfrag{(b)}{(b)}
\epsfig{height=20mm,file=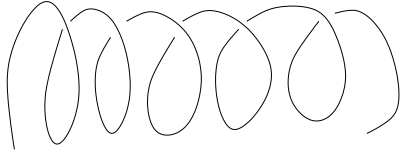}
\end{center}
\caption{A phone wire is locally convex in $\Ss^3$}
\label{fig:phonewire}
\end{figure}

For $n = 2$, $\xi$ is a circle;
for $n = 3$, $\xi$ turns around in one plane
while it turns around at a faster rate in another plane:
for suitable values of $a_i$ and $c_i$, $\xi$ looks
like a phone wire (see Figure \ref{fig:phonewire}).

\begin{proof}
A straight-forward calculation gives,
for $n+1 = 2k$,
\[ \det(\xi(t), \xi'(t), \ldots, \xi^{(n)}(t)) =
\left( \prod_i (c_i^2 a_i) \right)
\left( \prod_{i < j} \left((a_i - a_j)^2 (a_i + a_j)^2\right) \right) > 0 \]
and for $n+1 = 2k+1$,
\[ \det(\xi(t), \xi'(t), \ldots, \xi^{(n)}(t)) =
c_0 \left( \prod_i (c_i^2 a_i^3) \right)
\left( \prod_{i < j} \left((a_i - a_j)^2 (a_i + a_j)^2\right) \right) > 0. \]
Alternatively, we can compute $\Xi = \F_\xi: [0,1] \to SO_{n+1}$
and its logarithmic derivative $\Lambda: [0,1] \to so_{n+1}$:
it turns out to be a constant element of $\fT \subset so_{n+1}$.
In general, if $Q \in SO_{n+1}$ and $\gamma$ is locally convex
then so is $Q\gamma$: thus, $\xi_1$ is locally convex.

Conversely, assume $\tilde\xi \in \Little\Ss^n$
and that $\Lambda_{\tilde\xi}$ is constant equal to $B \in \fT$;
we then have $\tilde\Xi(t) = \F_{\tilde\xi}(t) = \exp(tB)$.
The eigenvalues of $B$ are all on the imaginary axis 
and therefore of the form $\pm a_j i$, $a_j > 0$,
plus a $0$ in case $n+1$ is odd.
Thus there exists $Q \in SO_{n+1}$ such that
\[ B = \begin{cases}
Q \diag\left(
\begin{pmatrix} 0 & -a_1 \\ a_1 & 0 \end{pmatrix}, 
\begin{pmatrix} 0 & -a_2 \\ a_2 & 0 \end{pmatrix}, \ldots,
\begin{pmatrix} 0 & -a_k \\ a_k & 0 \end{pmatrix} 
\right) Q^T,  & n = 2k-1, \\
Q \diag\left(0,
\begin{pmatrix} 0 & -a_1 \\ a_1 & 0 \end{pmatrix}, 
\begin{pmatrix} 0 & -a_2 \\ a_2 & 0 \end{pmatrix}, \ldots,
\begin{pmatrix} 0 & -a_k \\ a_k & 0 \end{pmatrix} 
\right) Q^T, & n = 2k. \end{cases}
\]
Write 
\[ X(s) = \begin{pmatrix} \cos(s) & -\sin(s) \\ \sin(s) & \cos(s) \end{pmatrix}
= \exp\left( s \begin{pmatrix} 0 & -1 \\ 1 & 0 \end{pmatrix} \right). \]
Thus, according to parity we have:
\[
\tilde\Xi(t) = \begin{cases}
Q \diag\left(X(a_1 t), X(a_2 t), \ldots, X(a_k t) \right) Q^T, \\
Q \diag\left(0, X(a_1 t), X(a_2 t), \ldots, X(a_k t) \right) Q^T, 
\end{cases} \]
and therefore
\[
\tilde\xi(t) =
Q \diag\left([0,] X(a_1 t), X(a_2 t), \ldots, X(a_k t) \right) v_0,
\quad v_0 = Q^T e_1. 
\]
Up to multiplication by a matrix of the form
$\diag\left( [1,] X(\theta_1), \ldots, X(\theta_k) \right)$,
$v_0$ can be assumed to be of the form
$v_0 = ([c_0,]c_1,0,\ldots,c_k,0)$ for $c_i \ge 0$
(with a corresponding change of $Q$).
The formulas in the previous paragraph of this proof
indicate that the parameters
$a_i$ and $c_i$ must be positive and that the $a_i$'s must be pairwise distinct
for $\tilde\xi$ to be locally convex, as desired.
The other claims are easy.
\end{proof}

The space of smooth curves is not the most convenient, however;
we use the above correspondence to define our favorite space of curves:
if we consider
$\Lambda \in L^2([0,1],\fT) \subset L^2([0,1],so_{n+1})$
we can solve the initial value problem $(\ast)$ and determine  
  $\Gamma: [0,1] \to SO_{n+1}$ and $\gamma(t) = \Gamma(t) e_1$.
Notice that the curve $\gamma$  constructed in this way from $\Lambda \in L^2$
belongs to $H^1([0,1],\RR^{n+1})$
but the concept of local convexity does not make sense
for all curves $\gamma: [0,1] \to \Ss^n$
with $\gamma \in H^1([0,1],\RR^{n+1})$.
A minor inconvenience is that $L^2([0,1],\fT)$ is not a Hilbert manifold;
we resolve this problem by defining a diffeomorphism
$\phi: \fT \to \RR^n$ with $j$-th coordinate
$\phi_j(T) = g(T_{j+1,j})$, $g(x) = x - 1/x$.
Given $\alpha \in L^2([0,1],\RR^n)$ we set 
$\Lambda = \phi^{-1} \circ \alpha$ and $\Gamma$ as above,
thus defining a space $\hat\Little\Ss^n$ of Jacobian curves
and an explicit diffeomorphism $L^2([0,1],\RR^n) \equiv \hat\Little\Ss^n$.
There are sometimes advantages in working with $\hat\Little\Ss^n(z)$
rather than $\Little\Ss^n(z)$: for instance, multiplication
(in $\Spin_{n+1}$) allows for a sort of superposition.
This will be useful later in the paper.

Given $Q \in SO_{n+1}$ let
$\hat\Little\Ss^n(z) \subset \hat\Little\Ss^n$ be the space
of Jacobian curves $\Gamma: [0,1] \to SO_{n+1}$
with $\Gamma(0) = I$, $\Gamma(1) = Q$.
Finally, we use the map $\tF$ to define the spaces
$\Little\Ss^n$ and $\Little\Ss^n(Q)$, which are now Hilbert manifolds.

Recall also that two Hilbert manifolds are diffeomorphic
if and only if they are homeomorphic which, in its turn, holds 
if and only if they are weakly homotopically equivalent (\cite{BH}).
Besides the Hilbert manifold structure,
the above definition of the spaces $\Little\Ss^n(Q)$ (and other spaces)
has the advantage  
of allowing discontinuities in $\Lambda_\gamma$
thus bypassing the need for a roundabout smoothening process.
The following result is now trivial.

\begin{lemma}
\label{lemma:Lcontract}
The space $\Little\Ss^n$ is contractible.
\end{lemma}

\bigbreak

\section{Bruhat cells and the Coxeter-Weyl group $B_{n+1}$}

As a first step, for any fixed  dimension $n$ 
we reduce  Problem~\ref{prob:main} 
to consideration of only finitely many different
values of $Q$ and $z$ using well-known group actions. 
The key observation here is that if $\gamma_1: [0,1] \to \Ss^n$
is positive locally convex and $A \in \RR^{n \times n}$ has positive
determinant then both $A\gamma_1: [0,1] \to \RR^{n+1}$
and $\gamma_2: [0,1] \to \Ss^n$
with $\gamma_2(t) = \widehat{A\gamma_1(t)} = A\gamma_1(t)/|A\gamma_1(t)|$
are positive locally convex.

Let $\U_{n+1}^+$ be the group of real upper-triangular
matrices with positive diagonal
and $\U_{n+1}^1 \subset \U_{n+1}^+$
be the subgroup of matrices with diagonal entries equal to one.
Consider the action of $\U_{n+1}^1$ on  $GL_{n+1}(\mathbb{R})$
by conjugation:
in what follows we will refer to the action
of $\U_{n+1}^1$ on different spaces as the \textit {Bruhat action}. 
This action induces the action of  $\U_{n+1}^1$ on
$SO_{n+1}$ as the postcomposition of the conjugation
with the orthogonalization.
In other words,  
$B(U,Q) = UQU'$ where $U'$ is the only matrix in $\U_{n+1}^+$
such that $UQU' \in SO_{n+1}$;
thus, $B(U,Q)$ is obtained from $UQ$ by Gram-Schmidt.
It is well-known that the Bruhat action on $SO_{n+1}$
has finitely many orbits.
These orbits are referred to as the \textit{Bruhat cells} of $SO_{n+1}$:
two orthogonal matrices $Q_1, Q_2 \in SO_{n+1}$ belong to the same
Bruhat cell if and only if there exist 
upper triangular matrices $U_1$, $U_2$ with positive
diagonal satisfying $U_1 Q_1 = Q_2 U_2$.
We denote the Bruhat cell of $Q \in SO_{n+1}$ by
$\Bruhat(Q) \subset SO_{n+1}$.

Let $B_{n+1} \subset O_{n+1}$ be the Coxeter-Weyl group
of signed permutation matrices and let
$B^{+}_{n+1} = B_{n+1} \cap SO_{n+1}$.
Let $\Diag^{+}_{n+1} \subset B^+_{n+1}$ be 
the subgroup of diagonal matrices with entries $\pm 1$ and determinant $1$;
thus $\Diag^{+}_{n+1}$ is isomorphic to $(\ZZ/2\ZZ)^{n}$.
Each Bruhat cell contains exactly one element $Q_0 \in B^{+}_{n+1}$
and is diffeomorphic to a cell
whose  dimension  equals the number of inversions of $Q_0$.
In other words, for each $Q \in SO_{n+1}$
there is a unique $Q_0 \in B^{+}_{n+1}$
such that there exist $U_1, U_2 \in \U_{n+1}^+$ satisfying $Q = U_1 Q_0 U_2$.

We recall an algorithm producing  $Q_0$ from a given $Q$; 
this algorithm will be used later,
particularly in the study of chopping.
Consider the first column of $Q$ and look for the lowest non-zero entry,
say $Q_{i1}$.
We first multiply $Q$ by a diagonal matrix $D \in \U_{n+1}^+$
to obtain a new matrix $\tilde Q = DQ$ for which $(\tilde Q)_{i1} = \pm 1$;
for simplicity, we may thus assume $Q_{i1} = \pm 1$.
We next perform row operations on $Q$ to clean the first column above row $i$:
in other words, we obtain $U_1 \in \U_{n+1}^+$ such that $\tilde Q = U_1 Q$
satisfies $\tilde Q e_1 = \pm e_i$;
again assume from now on that $Q e_1 = \pm e_i$.
Now perform column operations on $Q$ to clean row $i$
to the right of the first column, i.e.,
obtain $U_2 \in \U_{n+1}^+$ such that $\tilde Q = Q U_2$
satisfies $e_i^T \tilde Q = \pm e_1^T$.
Repeat the process for each column:
at the end of the process we obtain $Q_0 =  U_1 Q U_2$
($U_1, U_2 \in \U_{n+1}^+$)
for which there exists a permutation $\pi$ such that
$Q_0 e_i = \pm e_{\pi(i)}$.
In other words, $Q_0 \in B_{n+1}$;
since $\det Q_0  = \det U_1  \det Q  \det U_2  > 0$
we have $\det Q_0  = 1$ and $Q_0 \in B^{+}_{n+1}$.

Recall that $\Pi: \Spin_{n+1} \to SO_{n+1}$
is a group homomorphism and the double cover of $SO_{n+1}$
(for $n > 1$, this is the universal cover).
Let 
\[ \tilde B^{+}_{n+1} = \Pi^{-1}(B^{+}_{n+1}) \subset \Spin_{n+1}, \qquad
\widetilde\Diag^{+}_{n+1} = \Pi^{-1}(\Diag^{+}_{n+1})
\subset B^{+}_{n+1}; \]
the groups $\tilde B^{+}_{n+1}$ and $\widetilde\Diag^{+}_{n+1}$
are $\ZZ/2\ZZ$-central extensions of $B^{+}_{n+1}$ and $\Diag^{+}_{n+1}$,
respectively.
Notice that 
\[ |B^{+}_{n+1}| = 2^n(n+1)!, \qquad
|\tilde B^{+}_{n+1}| = 2^{n+1}(n+1)!, \qquad
|\widetilde\Diag^{+}_{n+1}| = 2^{n+1}; \]
for instance, $\widetilde\Diag^{+}_3$ is isomorphic to
the quaternion group $Q_8 = \{ \pm 1, \pm i, \pm j, \pm k \}$.

The Bruhat cell decomposition can be lifted to $\Spin_{n+1}$
where each cell contains a unique element of $\tilde B^{+}_{n+1}$.
Two elements of $SO_{n+1}$ or $\Spin_{n+1}$ 
will be called \textit{Bruhat equivalent}
if they belong to the same cell
in the corresponding Bruhat decomposition.
We will also write $\Bruhat(z) \subset \Spin_{n+1}$
for the Bruhat cell of $z \in \Spin_{n+1}$.

The Bruhat action of $\U_{n+1}^1$ on $SO_{n+1}$ induces the Bruhat action 
of $\U_{n+1}^1$ on the space $\Little\Ss^n$ as follows: 
given $\gamma \in \Little\Ss^n$ and $U \in \U_{n+1}^1$,
set $(B(U,\gamma))(t) = (B(U,\F_\gamma(t)))e_1$
(where $e_1 = (1,0,0,\ldots,0) \in \RR^{n+1}$).
Clearly, if $\gamma \in \Little\Ss^n(z)$
then $B(U,\gamma) \in \Little\Ss^n(B(U,z))$.
The following lemma is now easy.



\begin{lemma}
\label{lemma:bruhat}
If $Q_1, Q_2 \in SO_{n+1}$ (resp. $z_1, z_2 \in \Spin_{n+1}$)
are Bruhat equivalent then $\Little\Ss^n(Q_1)$ and $\Little\Ss^n(Q_2)$
(resp. $\Little\Ss^n(z_1)$ and $\Little\Ss^n(z_2)$) are homeomorphic.
\end{lemma}

This explicit homeomorphism will be used again
and we therefore introduce some notation.
Let $Q_1$ and $Q_2$ be as in the lemma:
there exists a matrix $U \in \U_{n+1}^1$
with $B(U,Q_1) = Q_2$ and therefore $B(U^{-1},Q_2) = Q_1$.
Define $\Bi_{Q_1,U,Q_2}: \Little\Ss^n(Q_1) \to \Little\Ss^n(Q_2)$
by $\Bi_{Q_1,U,Q_2}(\gamma) = B(U,\gamma)$
(for $\gamma \in \Little\Ss^n(Q_1)$).
Similarly define $\Bi_{z_1,U,z_2}: \Little\Ss^n(z_1) \to \Little\Ss^n(z_2)$.

\begin{proof} 
The map $\Bi_{Q_1,U,Q_2}$ is a homeomorphism
with inverse $\Bi_{Q_2,U^{-1},Q_1}$; the spin case is similar.
\end{proof}




\section{Time reversal}


In this and the two following sections
we introduce three natural operations
acting on $B_{n+1}^{+}$ and on the spaces of curves under consideration. 
 
The naive idea here would be to consider the curve
$t \mapsto \gamma(1-t)$;
this curve however may be negative locally convex
and has the wrong endpoints:
we show how to fix these minor problems.

Let $J_+ = \diag(1,-1,1,-1,\ldots) \in O_{n+1}$;
notice that $\det J_+  = (-1)^{n(n+1)/2}$.
For $Q \in SO_{n+1}$, define $\TR(Q) = J_+ Q^T J_+$.
The map $\TR: SO_{n+1} \to SO_{n+1}$ is an anti-automorphism 
which lifts to an anti-automorphism $\TR: \Spin_{n+1} \to \Spin_{n+1}$.
Indeed, given $z \in \Spin_{n+1}$ consider a path
$\tilde\alpha: [0,1] \to \Spin_{n+1}$
with $\tilde\alpha(0) = 1$, $\tilde\alpha(1) = z$;
let $\alpha = \Pi \circ \tilde\alpha$ and
$\beta: [0,1] \to SO_{n+1}$ with $\beta(t) = \TR(\alpha(t))$;
lift $\beta$ to define $\tilde\beta: [0,1] \to \Spin_{n+1}$
with $\tilde\beta(0) = 1$; define $\TR(z) = \tilde\beta(1)$.
The map is well defined:
two homotopic paths $\tilde\alpha_0$ and $\tilde\alpha_1$
yield homotopic paths $\alpha_0$ and $\alpha_1$;
the paths $\beta_0$ and $\beta_1$ are also homotopic
(apply $\TR$ to the homotopy) and therefore
$\tilde\beta_0(1) = \tilde\beta_1(1)$.

These two anti-automorphisms 
preserve the subgroups $\Diag^{+}_{n+1} \subset B^{+}_{n+1} \subset SO_{n+1}$
and $\widetilde\Diag^{+}_{n+1} \subset \tilde B^{+}_{n+1} \subset \Spin_{n+1}$.
In fact, the map $\TR: B^{+}_{n+1} \to B^{+}_{n+1}$
admits a simple combinatorial description:
the matrix $\TR(Q)$ is obtained from $Q$
by transposition and the change of sign of
all entries with $i+j$ odd.
We do not present a detailed combinatorial description of $\TR$
in $\tilde B^{+}_{n+1}$ but we record an observation for later use.

\begin{lemma}
\label{lemma:subtransit}
Let $s \in \ZZ$, $s \equiv n+1 \pmod 4$, $|s| < n+1$.
Then there exists $z \in \widetilde\Diag^{+}_{n+1}$
with $\trace(z) = s$ and $\TR(z) = - z$.
\end{lemma}

\begin{proof}
Let
\[ Q = \diag(-1,-1, \ldots, -1, -1, -1, 1, -1, 1, 1, \ldots, 1, 1) \]
and $z$ with $\Pi(z) = Q$. We claim that $z$ satisfies the claim;
in order to perform this computation we construct paths in $SO_{n+1}$
and lift them to $\Spin_{n+1}$.
Write
\[ X(t) = 
\begin{pmatrix} \cos(\pi t) & -\sin(\pi t) \\
\sin(\pi t) & \cos(\pi t) \end{pmatrix}, \qquad
Y(t) = \begin{pmatrix} \cos(\pi t) & 0 & -\sin(\pi t) \\ 0 & 1 & 0 \\
\sin(\pi t) & 0 & \cos(\pi t) \end{pmatrix}. \]
Take 
\[ \alpha(t) = \diag\left( X(t), \ldots, X(t), Y(t),
1, \ldots, 1\right), \]
with $((n+1-s)/4)-1$ small $X(t)$ blocks,
one large $Y(t)$ block followed by $((n+1+s)/2)$ ones.
Now lift the path $\alpha$ to $\tilde\alpha: [0,1] \to \Spin(n+1)$
with $\tilde\alpha(0) = \1$:
without loss of generality, $z = \tilde\alpha(1)$.
Clearly $\trace(z) = s$ and
\[ \TR\alpha(t) = \diag\left( X(t), \ldots, X(t), Y(-t),
1, \ldots, 1\right). \]
Thus $\alpha$ and $\TR\alpha$ are only different in two 
of the coordinates of the $Y$ block, where one makes a half-turn
one way and the other makes a half-turn the other way.
Thus $\TR(z) = -z$, as required.
\end{proof}

Notice that the map $\TR: so_{n+1} \to so_{n+1}$
given by $\TR(X) = J_+ X^T J_+$
satisfies $\TR(X) = X$ for $X \in \fT$.
For $\gamma \in \Little\Ss^n(Q)$,
define its  \textit{time reversal} by 
\[ \gamma^{\TR}(t) = J_+ Q^T \gamma(1-t) \]
where $Q^T$ is the transpose of $Q$.

\begin{lemma}\label{lemma:TR}
For any $\gamma \in \Little\Ss^n(Q)$ we have
$\gamma^{\TR} \in \Little\Ss^n(\TR(Q))$.
Furthermore,
\[ \F_{\gamma^{\TR}}(t) = J_+ Q^T \F_{\gamma}(1-t) J_+; 
\quad
\Lambda_{\gamma^{\TR}}(t) = \Lambda_\gamma(1-t). \]
In particular, if $\xi \in \Little\Ss^n(I)$ is a locally convex curve
for which $\Lambda_\xi$ is constant then $\xi_1^{\TR} = \xi_1$.

Time reversal yields explicit homeomorphisms
\[ \Little\Ss^n(Q) \approx \Little\Ss^n(\TR(Q)), \qquad
\Little\Ss^n(z) \approx \Little\Ss^n(\TR(z)). \]
\end{lemma} 

\begin{proof}
Consider a smooth curve $\gamma \in \Little\Ss^n(Q)$;
we must check that $\gamma^{\TR}$ is positive locally convex:
we have
\[ (\gamma^{\TR})^{(j)}(t) = (-1)^{j} J_+ Q^T \gamma^{(j)}(1-t) \]
and therefore
\begin{gather*}
\det\left((\gamma^{\TR})(t), (\gamma^{\TR})'(t), \ldots,
(\gamma^{\TR})^{(n)}(t)\right) = \\
= (-1)^{n(n+1)/2} \det J_+  \det Q^T  
\det\left(\gamma(1-t), \gamma'(1-t), \ldots,
\gamma^{(n)}(1-t)\right) =  \\
= \det\left(\gamma(1-t), \gamma'(1-t), \ldots,
\gamma^{(n)}(1-t)\right) > 0.
\end{gather*}
We must now check that $\F_{\gamma^{\TR}}(0) = I$ and
$\F_{\gamma^{\TR}}(1) = \TR(Q)$.
Recall that 
\begin{align*}
\left(\gamma(t), \gamma'(t), \ldots, \gamma^{(n)}(t)\right) 
&=  \F_{\gamma}(t)  R_0(t),\\
\left((\gamma^{\TR})(t), (\gamma^{\TR})'(t), \ldots,
(\gamma^{\TR})^{(n)}(t)\right)
&= \F_{\gamma^{\TR}}(t) R_1(t)
\end{align*}
where $R_0$ and $R_1$ are upper triangular matrices
with positive diagonal.
We have
\begin{gather*}
\left((\gamma^{\TR})(t), (\gamma^{\TR})'(t), \ldots,
(\gamma^{\TR})^{(n)}(t)\right) = \\
\\ = J_+ Q^T \left(\gamma(1-t), \gamma'(1-t), \ldots,
\gamma^{(n)}(1-t)\right) J_+
\end{gather*}
and therefore
\begin{gather*}
\F_{\gamma^{\TR}}(t) = 
J_+ Q^T \F_{\gamma}(1-t)  R_0(1-t) J_+ R_1^{-1}(t) = \\
= (J_+ Q^T \F_{\gamma}(1-t) J_+)(J_+ R_0(1-t) J_+ R_1^{-1}(t)).
\end{gather*}
Since $(J_+ Q^T \F_{\gamma}(1-t) J_+) \in SO_{n+1}$
and $(J_+ R_0(1-t) J_+ R_1^{-1}(t))$
is upper triangular with positive diagonal
we have
\[ \F_{\gamma^{\TR}}(t) = J_+ Q^T \F_{\gamma}(1-t) J_+; \]
in particular,
\begin{align*}
\F_{\gamma^{\TR}}(0) &= J_+ Q^T \F_{\gamma}(1) J_+ = I, \\
\F_{\gamma^{\TR}}(1) &= J_+ Q^T \F_{\gamma}(0) J_+ = \TR(Q).
\end{align*}
This completes the proof of the first claim
and of the first identity for smooth $\gamma$;
the second identity follows by taking derivatives
and the final claim is now easy.
The identities are extended to the general case
(i.e., $\gamma$ not necessarily smooth) by continuity,
thus completing the proof.
\end{proof}

\section{Arnold duality}

Let $A \in B^{+}_{n+1} \subset SO_{n+1}$
be the anti-diagonal matrix with entries
$(A)_{i,n+2-i} = (-1)^{(i+1)}$;
for instance, for $n = 2$ and $n = 3$ we have, respectively,
\[
A = \begin{pmatrix}
0 & 0 & 1 \\ 0 & -1 & 0 \\ 1 & 0 & 0 \end{pmatrix}, \qquad
A = \begin{pmatrix}
0 & 0 & 0 & 1 \\ 0 & 0 & -1 & 0 \\
0 & 1 & 0 & 0 \\ -1 & 0 & 0 & 0 \end{pmatrix}; \]
this matrix will appear in several places below.
Define an automorphism $\AD$ of $SO_{n+1}$ by $\AD(Q) = A^T QA$;
notice that the subgroup $B^{+}_{n+1} \subset SO_{n+1}$
is invariant under this automorphism.
As before, lift this automorphism to define an automorphism
(also called $\AD$) of $\Spin_{n+1}$ and $\tilde B^{+}_{n+1}$.
The combinatorial description of $\AD$ on $B^{+}_{n+1}$ is the following:
rotate $Q$ by a half-turn
(meaning that the $(i,j)$-th entry of $Q$
becomes the $(n-i+2,n-j+2)$-th entry of the new matrix) 
and change signs of all entries with $i+j$ odd.
Notice that the map $\AD: so_{n+1} \to so_{n+1}$
given by $\AD(X) = A^T X A$ takes $\fT$ to itself (as a set),
but reverts the order of the subdiagonal entries.

For $\gamma \in \Little\Ss^n(Q)$, define its  \textit{Arnold dual} as 
\[ \gamma^{\AD}(t) = \AD(\F_\gamma(t))e_1. \]
It turns out that this operation is just the usual projective duality
between oriented hyperplanes and unit vectors in disguise
(comp \cite {Arnold}).

\begin{lemma}\label{lemma:AD}
For any $\gamma \in \Little\Ss^n(Q)$ one has 
that $\gamma^{\AD} \in \Little\Ss^n(\AD(Q))$.
Furthermore,
\[ \F_{\gamma^{\AD}}(t) = \AD(\F_\gamma(t)), \quad
\Lambda_{\gamma^{\AD}}(t) = \AD(\Lambda_\gamma(t)). \]
Arnold duality gives explicit homeomorphisms
\[ \Little\Ss^n(Q) \approx \Little\Ss^n(\AD(Q)), \quad
\Little\Ss^n(z) \approx \Little\Ss^n(\AD(z)). \]
\end{lemma}

\begin{proof}
We must first check that if $\gamma \in \Little\Ss^n$ is smooth
then $\gamma^{\AD}$ as defined above also belongs to $\Little\Ss^n$.
Consider $\tilde\Gamma: [0,1] \to SO_{n+1}$ given by
\[ \tilde\Gamma(t) = \AD(\F_\gamma(t)) = A^T \F_\gamma(t) A. \]
We have
\[ (\tilde\Gamma(t))^{-1} {\tilde\Gamma}'(t) =
A^T (\F_\gamma(t))^{-1} A A^T \F_\gamma'(t) A =
A^T \Lambda_\gamma(t) A = \AD(\Lambda_\gamma(t)) \in \fT; \]
by Lemma \ref{lemma:Frenet},
$\tilde\Gamma = \F_{\tilde\gamma}$
for $\tilde\gamma \in \Little\Ss^n$:
thus $\gamma^{\AD} = \tilde\gamma \in \Little\Ss^n$,
completing our first check.
The formulas for $\F_{\gamma^{\AD}}$ and $\Lambda_{\gamma^{\AD}}$
have also been proved for smooth $\gamma$ and therefore,
by continuity, for all $\gamma \in \Little\Ss^n$.
The formulas imply that
if $\gamma \in \Little\Ss^n(Q)$ then $\gamma^{\AD} \in \Little\Ss^n(\AD(Q))$.
The final claim is now easy.
\end{proof}

\section{Chopping operation}

The first two operations corresponded to $\ZZ/2\ZZ$-symmetries
in $\Little\Ss^{n+1}$;
our third operation is quite different,
loosely corresponding to taking $\gamma \in \Little\Ss^{n+1}$
and chopping off a small tip at the end.
We again start with algebra and combinatorics.

For a signed permutation $Q \in B^{+}_{n+1}$
and a pair of indices $(i,j)$ with $(Q)_{(i,j)} \ne 0$
define $\NE(Q,i,j)$ to  be the number of pairs $(i',j')$ with
$i' < i$, $j' > j$ and $(Q)_{(i',j')} \ne 0$.
In other words,  $\NE(Q,i,j)$ is  
the number of nonzero entries of $Q$ in the northeast quadrant. 
Also set  
$$\SW(Q,i,j) = \NE(Q^T,j,i).$$ 
It is easy to check that for all $Q$ one has $\NE(Q,i,j) - \SW(Q,i,j) = i-j$.

Using the above notation define 
\[ \delta_i(Q) = (Q)_{(i,j)} (-1)^{\NE(Q,i,j)} \]
where $j$ is the only index for which $(Q)_{(i,j)} \ne 0$.
Additionally, define
\[ \Delta(Q) = \diag(\delta_1(Q), \delta_2(Q), \ldots, \delta_{n+1}(Q)),
\quad \text{and} \quad  \trd(Q) = \trace(\Delta(Q)). \]

\begin{lemma}\label{lm:5}
$\det Q  = \det(\Delta(Q))$.
\end{lemma}

\begin{proof}
Indeed, let $\pi$ be the permutation such that $\pi(j) = i$
if $j$ is the only index for which $(Q)_{(i,j)} \ne 0$.
Then
\begin{gather*}
\det(\Delta(Q)) = \prod_i \delta_i(Q)
= \left(\prod_i (Q)_{(i,j)}\right) (-1)^{\sum_i \NE(Q,i,j)} \\
= \left(\prod_i (Q)_{(i,j)}\right)
(-1)^{|\{ (i,i'); i' < i, \pi^{-1}(i') > \pi^{-1}(i) \}|}
= \det Q.
\end{gather*}
\end{proof} 

Thus $\Delta$ is a function from $B^{+}_{n+1}$
to $\Diag^{+}_{n+1} \subset B^{+}_{n+1}$.
Notice that $\Delta(Q) = Q$ for any $Q \in \Diag^{+}_{n+1}$.
We extend $\Delta$ to a function from $SO_{n+1}$ to $\Diag^{+}_{n+1}$ 
by declaring that if $Q$ and $Q'$ are Bruhat equivalent
then $\Delta(Q) = \Delta(Q')$; we similarly extend the function $\trd(Q)$ to $SO_{n+1}$.
The map $\Delta: SO_{n+1} \to \Diag^{+}_{n+1}$ is a projection
(in the sense that $\Delta(\Delta(Q)) = \Delta(Q)$)
and therefore defines a partition of $SO_{n+1}$
into $2^n$ classes of the form $\Delta^{-1}(Q)$, $Q \in \Diag^{+}_{n+1}$.
Furthermore, if $Q \in \Diag^{+}_{n+1}$,
we have $\Delta(QQ') = Q \Delta(Q')$ so that
a class $\Delta^{-1}(Q)$ is a fundamental domain
for the action of $\Diag^{+}_{n+1}$ on $SO_{n+1}$ by multiplication.

Let $A$ be the  matrix used in the definition of Arnold duality.
Notice that $\Delta(A) = I$
and therefore $\Delta(QA) = Q$ for all $Q \in \Diag^{+}_{n+1}$.
For $Q\in SO_{n+1}$, its \textit{chopping} is defined by $\chop(Q)=\Delta(Q)A$.
Thus the Bruhat equivalence class of $\chop(Q)$ is an open set,
dense in $\Delta^{-1}(\Delta(Q)) = \chop^{-1}(\chop(Q))$.
The maps $\Delta$ and $\chop$
as well as the functon $\trd: B^{+}_{n+1} \to \ZZ$
will  play a crucial role in our argument.   
(Notice that $\Delta$ is {\bf not} a group homomorphism).

Let us present a geometric interpretation for $\Delta$ and $\chop$. 
For $\gamma \in \Little\Ss^n(Q)$ and $\epsilon > 0$,
we define the \textit{naive chop} of $\gamma$ by $\epsilon$ as 
\[ \chop_\epsilon(\gamma)(t) = \gamma((1-\epsilon)t). \]
A straightforward computation gives 
\[ \F_{\chop_\epsilon(\gamma)}(t) = \F_\gamma((1-\epsilon)t), \quad
\Lambda_{\chop_\epsilon(\gamma)}(t) =
(1-\epsilon) \Lambda_\gamma((1-\epsilon)t); \]
in particular,
\[ \F_{\chop_\epsilon(\gamma)}(1) = \F_\gamma(1-\epsilon). \]
The inconvenience here is that if $\epsilon>0$ is fixed
and $\gamma$ varies over the whole $\Little\Ss^n(Q)$
we have no control of $\F_{\chop_\epsilon(\gamma)}(1)$,
the final frame of $\chop_\epsilon(\gamma)$.
The situation improves if we adapt the choice of $\epsilon$
depending on $\gamma$ and focus on Bruhat cells
instead of individual final frames.

\begin{lemma}\label{lm:chop} 
For any $Q \in SO(n+1)$ and for any $\gamma \in \Little\Ss^n(Q)$
there exists $\epsilon > 0$ such that 
for all $t \in (1-\epsilon,1)$
we have that $\F_\gamma(t) \in \Bruhat(\chop(Q))$.
\end{lemma}

In other words, given $\gamma \in \Little\Ss^n$
there exists $\tilde\epsilon > 0$ such that,
for all $\epsilon \in (0,\tilde\epsilon)$,
$\F_{\chop_\epsilon(\gamma)}(1)$ is Bruhat equivalent
to $\chop(\F_\gamma(1))$.

Before proving Lemma~\ref{lm:chop} we present
an illustrative  example for $n = 2$.
Take 
\[ Q = \begin{pmatrix} 0 & 1 & 0  \\ -1 & 0 & 0 \\ 0 & 0 & 1 \end{pmatrix}. \]
and  expand an arbitrary smooth curve $\gamma\in \Little\Ss^2(Q)$
in a Taylor series near $t = 1$. Using  $x = t-1$
we get  
$\gamma(x) \approx (x,-1,x^2/2)$ (up to higher order terms)  so that,
for $x \approx 0$,
\[ \F_\gamma(x) \approx
\begin{pmatrix} x & 1 & 0 \\ -1 & 0 & 0 \\ x^2/2 & x & 1  \end{pmatrix}. \]
We now apply the above algorithm 
to find $Q_0 \in B^{+}_{n+1}$ which is Bruhat equivalent 
to $\F_\gamma(x)$ when $x$ is a negative number with a small absolute value
(i.e., $Q_0 = U_1 \F_\gamma(x) U_2$).
We start at the $(3,1)$-th entry $x^2/2$, which is positive.
Thus, $(Q_0)_{3,1} = +1$.
We now concentrate on the $\SW$ (i.e., bottom left) $(2 \times 2)$-blocks
of $\F_\gamma(x)$ and $Q_0$:
since  $Q_0 = U_1 \F_\gamma(x) U_2$,
the signs of the determinants of these two blocks should be equal; 
since its original value equals  $-x > 0$,
the $(2,2)$-th entry of $Q_0$ equals $-1$.
Finally, the $(1,3)$-th entry must be set to $1$
for the whole determinant to be positive.
Summing up, if $\gamma \in \Little\Ss^2(Q)$
then there exists $\epsilon > 0$
such that for any $t \in (1-\epsilon,1)$
one has that $\F_\gamma(t)$ is Bruhat equivalent to
\[ \chop(Q) =
\begin{pmatrix} 0 & 0 & 1 \\ 0 & -1 & 0 \\ 1 & 0 & 0 \end{pmatrix}. \]
The general proof below follows the same idea.

\begin{proof}
As above, consider $Q \in B^{+}_{n+1}$
with associated permutation $\pi$,
so that $Q_{i,j} \ne 0$ if and only if $\pi(j) = i$.
Write a Taylor approximation $\F_\gamma(x) \approx M(x)$
where 
\[ (M(x))_{i,j} =
s_i g(\ell) x^\ell, \quad \pi^{-1}(i) = j+\ell, \quad (Q)_{i,j+\ell} = s_i, \]
where
\[ g(\ell) = \begin{cases}1/\ell!,&\ell \ge 0,\\0,&\ell < 0.\end{cases} \]
Let $M_k(x)$ be the $\SW$ $(k \times k)$-block of $M(x)$:
from the algorithm, we must show that, for small negative $x$, 
the matrix $M_k(x)$ is invertible and compute the sign of its determinant.

Write $M(x) = E G^\pi X^\pi(x) \tilde M \tilde X(x)$ for
\begin{gather*}
E = \diag(s_i), \quad
G^\pi = \diag(g(\pi^{-1}(i) - 1)), \quad
X^\pi(x) = \diag(x^{\pi^{-1}(i) - 1}), \\
(\tilde M)_{i,j} = (\pi^{-1}(i)-1)^{\underline{j-1}}, \quad
\tilde X(x) = \diag(x^{1-i}).
\end{gather*}
Here we use the notation $a^{\underline b} = a(a-1)\cdots(a-b+1)$.
Let $E_k, G^\pi_k, X^\pi_k(x)$ be the $\SE$ $k\times k$-blocks
of $E, G^\pi, X^\pi(x)$, respectively.
Similarly, let $\tilde M_k$ and $\tilde X_k(x)$ be the
$\SW$ and $\NW$ $k\times k$-blocks
of $\tilde M$ and $\tilde X(x)$, respectively.
We have $M_k(x) = E_k G^\pi_k X^\pi_k(x) \tilde M_k \tilde X_k(x)$
and therefore $\det M_k(x)$ is the product of the determinants
of these blocks. We must therefore determine the sign of the determinant
of each block.

For real numbers $a$ and $b$, we write $a \sim b$
if $a$ and $b$ have the same sign.
We have $\det E_k = \prod_{j \ge n-k+2} s_i$,
$\det G^\pi_k \sim 1$,
\begin{align*}
\det X^\pi_k(x) &= \prod_{j \ge n-k+2} (x^{\pi^{-1}(i) - 1}) =
x^{\left( \sum_{j \ge n-k+1} (\pi^{-1}(i) - 1) \right)} \\
&\sim (-1)^{\left(k+ \sum_{j \ge n-k+2} \pi^{-1}(i) \right)} 
\end{align*}
and $\det \tilde X_k(x) = x^{-k(k-1)/2} \sim (-1)^{k(k-1)/2}$.
In order to compute $\det \tilde M_k$, consider the Vandermonde matrix
$V^\pi$ with $(V^\pi)_{i,j} = (\pi^{-1}(i)-1)^{j-1}$;
notice that there exists $U \in \U_{n+1}^1$ with $V^\pi = \tilde M U$.
Let $V^\pi_k$ be the $\SW$ $k\times k$-block of $V^\pi$,
also a Vandermonde matrix. We have
\[ \det \tilde M_k = \det V^\pi_k =
\prod_{n-k+2 \le j < j' \le n+1} (\pi^{-1}(j') - \pi^{-1}(j)). \]
At this point we know that $\det M_k \ne 0$
(with the same sign for all small negative $x$)
and therefore there exists a diagonal matrix $\hat\Delta(Q) \in B^+_{n+1}$
such that $M(x)$ and $\hat M = \hat\Delta(Q) A$ are Bruhat equivalent.
Write $\hat\Delta(Q) = \diag(\hat\delta_i(Q))$;
we must compute $\hat\delta_i(Q)$.
Let $\hat M_k$ be the $SW$ $k \times k$-block of $\hat M_k$:
by Bruhat equivalence we have $\det \hat M_k \sim \det M_k(x)$;
by construction we have
$\det \hat M_k = (-1)^{kn} \prod_{j \ge n-k+2} \hat\delta_j(Q)$.
Thus $\hat\delta_{n-k+2}(Q) \sim (-1)^n \det M_k(x) \det M_{k-1}(x)$.

We have
\begin{align*}
\det E_k \det E_{k-1} &= s_{n-k+2} = Q_{n-k+2,\pi^{-1}(n-k+2)}, \\
\det X^\pi_k(x) \det X^\pi_{k-1}(x) &\sim (-1)^{\pi^{-1}(n-k+2) - 1}, \\
\det \tilde X_k \det \tilde X_{k-1} &\sim (-1)^{k-1}, \\
\det \tilde M_k \det \tilde M_{k-1} &\sim
\prod_{n-k+2 < j' \le n+1} (\pi^{-1}(j') - \pi^{-1}(n-k+2)) \\
&\sim (-1)^{\SW(Q,n-k+2,\pi^{-1}(n-k+2))}
\end{align*}
and therefore
\[ \hat\delta_{n-k+2}(Q) \sim (-1)^n Q_{n-k+2,\pi^{-1}(n-k+2)} 
(-1)^{\SW(Q,n-k+2,\pi^{-1}(n-k+2)) + k + \pi^{-1}(n-k+2)}. \]
Since both sides have absolute value $1$ the latter relation is actually an equality;
for $i = n-k+2$ and $j = \pi^{-1}(n-k+2)$ we then have
\[ \hat\delta_{i}(Q) =
Q_{i,j} (-1)^{\SW(Q,i,j) + i + j} =
Q_{i,j} (-1)^{\NE(Q,i,j)} = \delta_i(Q). \]
Thus $\hat\Delta(Q) = \Delta(Q)$ and we are done.
\end{proof} 

A geometric description of the situation is now more clear.  
The Bruhat cells of the form $\Bruhat(DA)$, $D \in \Diag^+_{n+1}$,
are disjoint open sets and their union is dense in $SO_{n+1}$.
The complement of this union is the disjoint union
of Bruhat cells of lower dimension.
Let $\Gamma: (-\epsilon,\epsilon) \to SO_{n+1}$ be a smooth Jacobian curve
(i.e., with $\Lambda(t) = (\Gamma(t))^{-1} \Gamma'(t) \in \fT$
for all $t \in (-\epsilon,\epsilon)$):
if $\Gamma(0)$ does not belong to a top-dimensional Bruhat cell
then the function $\chop$ and Lemma \ref{lm:chop}
tell us in which cell $\Gamma(t)$ falls for $t < 0$, $|t|$ small.
In other words, provided you follow a Jacobian curve you can only arrive at a given low-dimensional
Bruhat cell from one of the adjacent top-dimensional cells.

As discussed above, the decomposition into Bruhat cells
lifts of $\Spin_{n+1}$.
The above geometric characterization of $\chop$
thus also lifts to a map
$\chop: \Spin_{n+1} \to \tilde B^{+}_{n+1}$.
Let $\ba = \chop(\1)$ (so that $\Pi(\ba) = A$)
and define $\Delta: \Spin(n+1) \to \widetilde\Diag^{+}_{n+1}$
by $\chop(z) = \Delta(z) \ba$.   
We shall not attempt to give a combinatorial description
of $\Delta$ or $\chop$ in the spin groups.

We present yet another interpretation of the chopping operation. 
Let $\Gamma: (t_0 - c, t_0 + c) \to \Spin_{n+1}$
be a Jacobian curve.
Notice that if $\Gamma$ is Jacobian and $z \in \Spin_{n+1}$
then so is $z \Gamma$ (their logarithmic derivatives are equal).
Thus, Lemma \ref{lm:chop} can be extended to show
that for all $z \in \Spin_{n+1}$  one has that 
$z \Gamma(t_0 - \epsilon) \in \Bruhat(\chop(z \Gamma(t_0)))$
or
$\Gamma(t_0 - \epsilon) \in z^{-1} \Bruhat(\chop(z \Gamma(t_0)))$.
In particular, taking $z = (\Gamma(t_0))^{-1}$, we have
$\Gamma(t_0 - \epsilon) \in \Gamma(t_0) \Bruhat(\ba)$.
Conversely, given $Q_1 \in SO_{n+1}$ and
$Q_0 \in Q_1 \Bruhat(A)$ there exists
a globally Jacobian curve $\Gamma: [0,1] \to SO_{n+1}$
with $\Gamma(0) = Q_0$, $\Gamma(1) = Q_1$ 
(so that $\gamma: [0,1] \to \Ss^n$, $\gamma(t) = \Gamma(t) e_1$,
is globally convex).
The following statement thus follows from Lemma \ref{lm:chop}.

\begin{coro}
\label{coro:chopbruhat}
Given $Q \in SO_{n+1}$ there exists an open set $U \subset SO_{n+1}$
with $Q \in U$ and $U \cap (Q \Bruhat(A)) \subseteq \Bruhat(\chop(Q))$.
Similarly, given $z \in \Spin_{n+1}$
there exists an open set $U \subset \Spin_{n+1}$
with $z \in U$ and $U \cap (z \Bruhat(\ba)) \subseteq \Bruhat(\chop(z))$.
\end{coro}

\begin{proof}
The $SO_{n+1}$ case follows from the remarks above
together with Lemma \ref{lm:chop};
the $\Spin_{n+1}$ case is similar.
\end{proof}

The next statement is crucial in our consideration.

\begin{prop}
\label{prop:chop}
For any $z \in \tilde B^{+}_{n+1}$
there are homeomorphisms
$$\Little\Ss^{n+1}(z)
\approx \Little\Ss^{n+1}(\chop(z))
\approx \Little\Ss^{n+1}(\Delta(z)).$$
\end{prop}

We need a few preliminary constructions and results.
Consider a Jacobian curve $\Gamma_0: [0,1] \to SO_{n+1}$
with $\Gamma_0(0) = Q_0$ and $\Gamma_0(1) = Q_1$.  Define
$\Ci_{Q_0,\Gamma_0,Q_1}: \Little\Ss^{n+1}(Q_0) \to \Little\Ss^{n+1}(Q_1)$
by 
\[ (\Ci_{Q_0,\Gamma_0,Q_1}(\gamma))(t) = \begin{cases}
\gamma(2t),&t \le 1/2, \\ \Gamma_0(2t-1) e_1,&t \ge 1/2.\end{cases} \]

\begin{lemma}
\label{lemma:BiCi}
Consider a globally Jacobian curve $\Gamma_0: [0,1] \to SO_{n+1}$
whose image is contained in a Bruhat cell.
Let $Q_0 = \Gamma_0(0)$, $Q_1 = \Gamma_0(1)$ and
$U \in \U_{n+1}^1$ with $B(U,Q_0) = Q_1$.
Then the maps $\Bi_{Q_0,U,Q_1}$ and $\Ci_{Q_0,\Gamma_0,Q_1}$
are homotopic. In particular, $\Ci_{Q_0,\Gamma_0,Q_1}$
is a homotopy equivalence.
\end{lemma}

\begin{proof} 
Since $Q_0$ and $\Gamma(s)$ 
are in the same Bruhat cell  for all $s\in [0,1]$ we can define
a continuous function $\bU: [0,1] \to \U_{n+1}^1$
with $B(\bU(s),Q_0) = \Gamma(s)$, $\bU(0) = I$, $\bU(1) = U$.
Define $H: \Little\Ss^n(Q_0) \times [0,1] \to \Little\Ss^n(Q_1)$
by 
\[ (H(\gamma,s))(t) = \begin{cases}
\Bi_{Q_0,\bU(s),\Gamma(s)}(2t/(1+s)),&t \le (1+s)/2,\\
\Gamma_0(2t-1) e_1,&t \ge (1+s)/2. \end{cases} \]
The map $H$ produces  the desired homotopy from 
$\Ci_{Q_0,\Gamma_0,Q_1}$ to $\Bi_{Q_0,U,Q_1}$.
\end{proof} 

To prove Proposition~\ref{prop:chop}
we will also use the following previously known facts.

\begin{fact}[comp. Lemma $5$ in \cite{ShapiroM}] 
\label{fact:contractible}
For any $z \in \Spin_{n+1}$ the space $\Little\Ss^n(z)$
has two connected components if and only if 
there exists a globally convex curve in $\Little\Ss^n(z)$. 
One of these connected components
is the set of all globally convex curves in $\Little\Ss^n(z)$
and this connected component is contractible.
If $\Little\Ss^n(z)$ contains no globally convex curves  
then it is connected. 
\end{fact}

\begin{fact}[comp. Theorem 0.1 in \cite{BH}]
\label{fact:BH}
Let $M$ and $N$ be two topological Hilbert manifolds. Then 
any weak homotopy equivalence $f_0: N \to M$ is homotopic
to a homeomorphism $f_1: N \to M$.
\end{fact}


Let $z_1 \in \tilde B^+_{n+1}$ and
consider a smooth Jacobian curve
$\Gamma_{\aux}: [-\epsilon,\epsilon] \to \Spin_{n+1}$
with $\Gamma_{\aux}(0) = z_1$.
Choose $\epsilon$ sufficiently small so that the image of
$\Gamma_\aux([-\epsilon,0)) \subset
\Bruhat(\chop(z_1)) \cap (z_1 \Bruhat(\ba))$.
Let $z_0 = \Gamma_{\aux}(-\epsilon)$,
$\Gamma_0(t) = \Gamma_{\aux}(\epsilon(t-1))$.
Proposition \ref{prop:chop} now follows 
directly from the next lemma.

\begin{lemma}
\label{lemma:Cichop}
The map $\Ci_{z_0,\Gamma_0,z_1}:
\Little\Ss^n(z_0) \to \Little\Ss^n(z_1)$
is a weak homotopy equivalence.
\end{lemma}



\begin{proof} 
For $k$ a non-negative integer,
let $\alpha: \Ss^k \to \Little\Ss^{n+1}(z_1)$:
we construct $\tilde\alpha: \Ss^k \to  \Little\Ss^{n+1}(z_0)$
and a homotopy $H: \Ss^k \times [0,1] \to \Little\Ss^{n+1}(z_1)$
with $H(s,0) = \alpha(s)$, $H(s,1) = \Ci_{z_0,\Gamma_0,z_1}(\tilde\alpha(s))$.
By compactness and continuity, there exists $\epsilon_1 > 0$ such that
for all $s \in \Ss^k$ and for all $t \in [1-\epsilon_1,1)$
we have $\F_{\alpha(s)}(t) \in \Bruhat(\chop(z_1)) \cap (z_1 \Bruhat(\ba))$.
Again by compactness and continuity, there exists $\epsilon_2 > 0$,
$\epsilon_2 < \epsilon_1/2$, such that for all $s \in \Ss^k$
we have $\F_{\alpha(s)}(1-\epsilon_1) \in
\Bruhat(\chop(z_1)) \cap (\Gamma_0(1-\epsilon_2) \Bruhat(\ba))$.
Thus, for each $s \in \Ss^k$,
the space $X_s$ of globally Jacobian curves $\Gamma_s: [1-\epsilon_1,1]$
for which $\Gamma_s(1-\epsilon) = \F_{\alpha(s)}(1-\epsilon_1)$ and
$\Gamma_s(1) = z_1$ is non-empty
(since $\alpha(s)|_{[1-\epsilon_1,1]} \in X_s$) and therefore,
by Fact \ref{fact:contractible}, a contractible space.
Consider the subspace $Y_s \subset X_s$ of curves
for which $\Gamma_s(t) = \Gamma_0(t)$ for $t \ge 1-\epsilon_2$;
the condition 
$\F_{\alpha(s)}(1-\epsilon_1) \in (\Gamma_0(1-\epsilon_2) \Bruhat(\ba))$
implies that $Y_s$ is non-empty and Fact \ref{fact:contractible}
implies that is $Y_s$ also contractible.
We may therefore construct
a homotopy $H_1: \Ss^k \times [0,1] \to \Little\Ss^{n+1}(z_1)$
with $H(s,0) = \alpha(s)$, $H(s,\tilde s) \in X_s$
and $H(s,1) \in Y_s$.
In other words, we may assume without loss of generality
that there exists $\epsilon_2 > 0$ such that
$\alpha(s)(t) = \gamma_0(t)$ for all $s \in \Ss^k$ and $t > 1-\epsilon_2$.

Set $z_2 = \Gamma_0(1-\epsilon_2)$ and $\Gamma_2: [0,1] \to \Spin_{n+1}$
with $\Gamma_2(t) = \Gamma_0((1-\epsilon_2) + \epsilon_2 t)$. 
We may reparameterize the curves so that, for all $s$,
$\alpha(s)(1/2) = z_2$ and $\alpha(s)(t) = \Gamma_2(2t-1)$ for $t \ge 1/2$.
In other words, we may assume that
$\alpha(s) = \Ci_{z_2,\Gamma_2,z_1} \hat\alpha(s)$.
Set $\Gamma_3(t) = \Gamma_0(t/(1-\epsilon_2))$;
Lemma \ref{lemma:BiCi} tells us that
$\Ci_{z_0,\Gamma_3,z_2}: \Little\Ss^n(z_0) \to \Little\Ss^n(z_2)$
is a homotopy equivalence:
$\hat\alpha$ is therefore homotopic to 
$\Ci_{z_0,\Gamma_3,z_2} \circ \tilde\alpha$ 
for some $\tilde\alpha: \Ss^k \to \Little\Ss^n(z_0)$,
implying that $\alpha$ is homotopic to
$\Ci_{z_0,\Gamma_0,z_1} \circ \tilde\alpha$, as desired.
This completes the proof that
$\pi_k(\Ci_{z_0,\Gamma_0,z_1}):
\pi_k(\Little\Ss^n(z_0)) \to \pi_k(\Little\Ss^n(z_1))$
is surjective.

The proof that this map is injective is similar.
Let $\tilde\alpha: \Ss^k \to \Little\Ss^n(z_0)$ and
$\alpha = \Ci_{z_0,\Gamma_0,z_1} \circ \tilde\alpha$;
assume that $\tilde\alpha$ is homotopically trivial,
i.e., that there exists $H: \BB^{k+1} \to  \Little\Ss^n(z_1)$
with $H|_{\Ss^k} = \alpha$:
we need to prove that $\alpha$ is homotopically trivial.
As above, change $H$ so that $H(s)$ agrees with $\Gamma_0$
near $1$, i.e., we may assume $H$ to be of the form
$H = \Ci_{z_2,\Gamma_2,z_1} \circ \hat H$.
We therefore have that
$\Ci_{z_0,\Gamma_3,z_2} \circ \tilde\alpha$ 
is homotopically trivial.
Since 
$\Ci_{z_0,\Gamma_3,z_2}: \Little\Ss^n(z_0) \to \Little\Ss^n(z_2)$
is a homotopy equivalence we are done.
\end{proof}

\section{Proof of Theorem \ref{theo:main}}

First we reformulate  Theorem \ref{theo:main}
using the language of the prevous sections.

\begin{theo}
\label{theo:mainpp}
Let $Q_0, Q_1 \in SO_{n+1}$: if $\trd(Q_0) = \trd(Q_1)$ then
$\Little\Ss^n(Q_0)$ and $\Little\Ss^n(Q_1)$ are homeomorphic.

Let $z_0, z_1 \in \Spin_{n+1}$: if $\trd(z_0) = \trd(z_1)$ and
$|\trd(z_0)| \ne n+1$ then
$\Little\Ss^n(z_0)$ and $\Little\Ss^n(z_1)$ are homeomorphic.
\end{theo}

Theorem \ref{theo:main} follows directly from Theorem \ref{theo:mainpp}.
The condition $|\trd(z_0)| \ne n+1$ in the spin part is necessary:
for $n = 2$ and $\1, -\1 \in \Spin_3$ the two central elements
the spaces $\Little\Ss^2(\1)$ and $\Little\Ss^2(-\1)$ are not homeomorphic
since they have different numbers of connected components.

Recall that from Lemmas \ref{lemma:bruhat} and \ref{lemma:Cichop}
and Proposition \ref{prop:chop}
we already know that if $\Delta(Q_0) = \Delta(Q_1)$
then $\Little\Ss^n(Q_0)$ and $\Little\Ss^n(Q_1)$ 
(as well as  $\Little\Ss^n(\Delta(Q_0))$) are homeomorphic;
we have a similar result for the spin group.
We are therefore left to consider  the spaces $\Little\Ss^n(Q)$,
$Q \in \Diag^+_{n+1}$, and their spin counterparts.
A number of additional statements  are required for 
the proof of Theorem \ref{theo:mainpp}.

\begin{lemma}
\label{lemma:transit}
Let $D_0, D_1 \in \Diag^{+}_{n+1}$ with $\trd(D_0) = \trd(D_1)$.
Then there exists $Q \in B^{+}_{n+1}$ with
$\Delta(Q) = D_0$, $\Delta(\TR(Q)) = D_1$.
Thus $\Little\Ss^n(D_0)$ and $\Little\Ss^n(D_1)$ are homeomorphic. 
\end{lemma}

{\nobf Proof: }
Let $\pi$ be a permutation of $\{1,2,\ldots,n+1\}$
with $(D_1)_{\pi(i),\pi(i)} = (D_0)_{i,i}$ for all $i$.
Let $P$ be a permutation matrix with $(P)_{(i,j)} = 1$
if and only if $j = \pi(i)$.
Set $Q = D_0 \Delta(P) P$: 
we have $\Delta(Q) = D_0 \Delta(P) \Delta(P) = D_0$.
On the other hand, if $\pi(i) = j$,
we have $\delta_j(\TR(Q)) = \delta_i(Q)$ 
(from the proof of Lemma \ref{prop:adtr})
and therefore
$\delta_j(\TR(Q)) = (D_0)_{(i,i)} (\Delta(P))_{(i,i)} \delta_i(P)
= (D_0)_{(i,i)} = (D_1)_{j,j}$  and $\Delta(\TR(Q)) = D_1$. 
The last claim follows from Proposition \ref{prop:chop}. 
\qed

This completes the proof of Theorem \ref{theo:mainpp}
for the $SO_{n+1}$ case:
one judicious use of the equivalences
proved in the previous section is enough.
The spin case is slightly subtler:
it turns out that a single instance of the equivalences
is not enough,  which can be readily checked
by an exhaustive search in the case $n = 2$.
A small chain of consecutive instances of the equivalences are
therefore used.

\begin{lemma}
\label{lemma:tildetransit}
Let $z_0, z_1 \in \widetilde\Diag^{+}_{n+1}$
with $\trd(z_0) = \trd(z_1) \ne \pm(n+1)$.
Then there exist $w_0, w_1 \in \tilde B^{+}_{n+1}$ with
$\Delta(w_0) = z_0$, $\Delta(\TR(w_1)) = z_1$
and either $\Delta(\TR(w_0)) = \Delta(w_1)$
or $\Delta(\TR(w_0)) = \TR(\Delta(w_1))$.
Thus $\Little\Ss^n(z_0)$ and $\Little\Ss^n(z_1)$ are homeomorphic.
\end{lemma}

\begin{proof}
Take $s = \trd(z_0)$ and apply Lemma \ref{lemma:subtransit}
to obtain $z \in \widetilde\Diag^{+}_{n+1}$
with $\TR(z) = -z$.
Let $Q_0 = \Pi(z_0)$, $Q_1 = \Pi(z_1)$, $Q = \Pi(z)$.
By Lemma \ref{lemma:transit} there exist $P_0, P_1 \in B^{+}_{n+1}$
with $\Delta(P_0) = Q_0$, $\Delta(\TR(P_0)) = Q$,
$\Delta(P_1) = Q$, $\Delta(\TR(P_1)) = Q_1$.
Take $w_0, w_1 \in \tilde B^{+}_{n+1}$ with
$\Pi(w_0) = P_0$, $\Pi(w_1) = P_1$,
$\Delta(w_0) = z_0$, $\Delta(\TR(w_1)) = z_1$.
We have $\Delta(\TR(w_0)) = \pm z$ and $\Delta(w_1) = \pm z$ and 
we are done.
\end{proof}

Theorem \ref{theo:mainpp} follows directly from 
Lemmas \ref{lemma:transit} and \ref{lemma:tildetransit}.

It is natural to ask whether Theorem \ref{theo:mainpp} is the strongest
possible such statement, i.e., if spaces which it does not declare
homeomorphic are actually not homeomorphic.
We do not know the answer to this question
(see Problem \ref{prob:mainppisstrong} below)
but the following proposition shows that it is the strongest
result which follows (or follows directly)
from the remarks of the previous sections.

\begin{prop}
\label{prop:adtr}
For all $Q \in B^{+}_{n+1}$ we have $\trd(\AD(Q)) = \trd(\TR(Q)) = \trd(Q)$.
\end{prop}

{\nobf Proof: }
Assume $(Q)_{(i,j)} \ne 0$.
We have
\begin{gather*}
\delta_{n+2-i}(\AD(Q)) =
(\AD(Q))_{(n+2-i,n+2-j)} (-1)^{\NE(\AD(Q),n+2-i,n+2-j)}\\
= (-1)^{(i+j)} Q_{(i,j)} (-1)^{\SW(Q,i,j)} =
Q_{(i,j)} (-1)^{\NE(Q,i,j)} = \delta_i(Q)
\end{gather*}
and
\begin{gather*}
\delta_j(\TR(Q)) =
(\TR(Q))_{(j,i)} (-1)^{\NE(\TR(Q),j,i)} \\
= (-1)^{(i+j)} (Q)_{(i,j)} (-1)^{\SW(Q,i,j)} =
Q_{(i,j)} (-1)^{\NE(Q,i,j)} = \delta_i(Q).
\end{gather*}
The proposition now follows. 
\qed

\section{Proof of Theorem \ref{theo:dois}}

Our nearest goal is  to  prove Theorem~\ref{theo:dois}(i),
i.e. the fact  that the inclusion
$\hat\Little\Ss^n(z) \subset \Omega\Spin_{n+1}(z)$ is homotopically surjective
for all $z$ and then to settle Theorem~\ref{theo:dois}(ii),
i.e. that this inclusion is  a homotopy equivalence if $\Pi(z) = \pm J_+$.

Recall that the group $SO_{n+1} \subset \RR^{(n+1)\times(n+1)}$
has a natural Riemann metric and $\Spin_{n+1}$ inherits it via $\Pi$.
With this metric,
let $r_{n+1} > 0$ be the injectivity radius of the exponential map,
i.e., $r_{n+1}$ is such that if $z_0,z_1 \in \Spin_{n+1}$,
$d(z_1,z_2) < r_{n+1}$, then there exists a unique shortest
geodesic $g_{z_0,z_1}: [0,1] \to \Spin_{n+1}$
(parametrized by a constant multiple of arc length)
joining $z_0$ and $z_1$ so that $g_{z_0,z_1}(i) = z_i$, $i = 0, 1$.

We will need another technical lemma.

\begin{lemma}
\label{lemma:fast}
Let $K$ be a smooth compact manifold and
$\alpha: K \times [0,1] \to \Spin_{n+1}$ be a smooth function
and write $\alpha_s(t) = \alpha(s,t)$.
Then there exists $\xi_\fast \in \Little\Ss^n(\1)$ 
and corresponding $\Xi_\fast \in \hat\Little\Ss^n$ such that
$\xi_\fast^{\TR} = \xi_\fast$ and
the curves $\gamma_s(t) =  \alpha_s(t) \xi_\fast(t)$
are positive locally convex for all $s \in K$.
Furthermore, given $\epsilon > 0$, $\epsilon < r_{n+1}$,
we may assume that
\[ d(\tF_{\gamma_s}(t), \alpha_s(t) \Xi_\fast(t)) < \epsilon \]
for all $s \in K$, $t \in [0,1]$.
\end{lemma}

\begin{figure}[ht]
\begin{center}
\epsfig{height=30mm,file=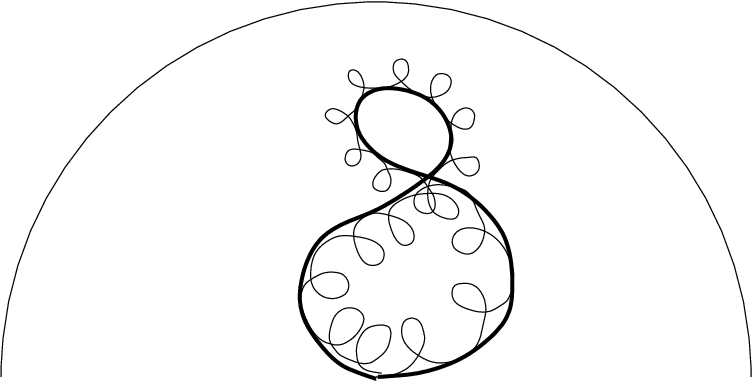}
\end{center}
\caption{Approximating a curve by a locally convex curve}
\label{fig:addloop}
\end{figure}

The intuitive picture here, at least for $n = 2$,
is that an arbitrary curve $\gamma: \Ss^1 \to \Ss^n$
can be replaced by a phone wire,
a locally convex curve which in some sense follows $\gamma$
while quickly rotating in a transversal direction  to guarantee local convexity
(see Figure \ref{fig:addloop}).

{\nobf Proof: }
Take $\xi_1 \in \Little\Ss^n(\1)$ as in Lemma \ref{lemma:phonewire}.
We claim that $\xi_\fast(t) = \xi_1(Nt)$
satisfies the lemma for a sufficiently large integer $N$.
Notice that $\xi^{(k)}_\fast(t) = N^k \xi_1^{(k)}(Nt)$
and 
\begin{align*}
\gamma_s^{(k)}(t) &= N^k \left(
\alpha_s(t) \xi_1^{(k)}(Nt) + \cdots +
\frac{1}{N^j} \binom{k}{j} \alpha_s^{(j)}(t) \xi_1^{(k-j)}(Nt) + \cdots 
\right) \\
&= N^k ( \alpha_s(t) \xi_1^{(k)}(Nt) + E_k(N,s,t) )
\end{align*}
where $E_k(N,s,t)$ tends to $0$ when $N\to \infty$.
Since
\[ \det(\alpha_s(t) \xi_1(Nt), \ldots, 
\alpha_s(t) \xi_1^{(n)}(Nt)) =
\det(\xi_1(Nt), \ldots, \xi_1^{(n)}(Nt)) \]
is positive and bounded away from $0$ it follows that $\gamma_s$
is indeed locally convex for sufficiently large $N$.
Furthermore, the identities
\begin{gather*}
\begin{pmatrix} \alpha_s(t) \xi_1(Nt) + E_0(N,s,t) &
 \cdots & \alpha_s(t) \xi_1^{(n)}(Nt) + E_n(N,s,t) \end{pmatrix}
= \F_{\gamma_s}(t) R_{\gamma_s}(t) \\
\begin{pmatrix} \alpha_s(t) \xi_1(Nt) &
 \cdots & \alpha_s(t) \xi_1^{(n)}(Nt) \end{pmatrix}
= \alpha_s(t) \Xi_\fast(t) R_{\xi_\fast}(t),
\end{gather*}
where $R_{\gamma_s}(t)$ and $R_{\xi_\fast}(t)$
are upper triangular matrices with positive diagonals,
show that $d(\tF_{\gamma_s}(t), \alpha_s(t) \Xi_\fast(t))$
can be made arbitrarily small by choosing large $N$.
\qed

\begin{prop}
\label{prop:surjective}
For any $z \in \Spin_{n+1}$ 
the inclusion $\hat\Little\Ss^n(z) \subset \Omega\Spin_{n+1}(z)$
is homotopically surjective.
In other words, given $\alpha_0: \Ss^k \to \Omega\Spin_{n+1}(z)$
there exists a homotopy in $\Omega\Spin_{n+1}(z)$ from $\alpha_0$ to
$\alpha_1: \Ss^k \to \hat\Little\Ss^n(z) \subset \Omega\Spin_{n+1}(z)$.
\end{prop}

\begin{proof}
Write $\alpha_0(s;t) = \alpha_0(s)(t)$.
Assume without loss of generality that $\alpha_0$ is smooth
if interpreted as $\alpha_0: \Ss^k \times [0,1] \to \Spin_{n+1}$.
We may also assume that $\alpha_0$ is flat at both $t=0$ and $t=1$, i.e.,
that $(\alpha_0(s))^{(m)}(t) = 0$ for $t \in \{0,1\}$,
for all $s \in \Ss^k$ and all $m>0$.
By Lemma \ref{lemma:fast},
there exist $\xi_\fast \in \Little\Ss^n(\1)$
and corresponding $\Xi_\fast \in \hat\Little\Ss^n(\1)$
such that 
\[ d(\tF_{\gamma_s}(t), \alpha_0(s;t) \Xi_\fast(t))
< \frac{r_{n+1}}{2} \]
for all $s \in \Ss^k$, $t \in [0,1]$;
here, as in  Lemma \ref{lemma:fast},
$\gamma_s(t) =  \alpha(s;t) \xi_\fast(t)$.
The flatness condition guarantees that
\[ \tF_{\gamma_s}(0) = \alpha_0(s;0) \Xi_\fast(0) = \1, \qquad
\tF_{\gamma_s}(1) = \alpha_0(s;1) \Xi_\fast(1) = z. \]
Recall that $\Xi_\fast \in \Omega\Spin_{n+1}(\1)$:
let $H:[0,1] \to \Omega\Spin_{n+1}(\1)$ be a homotopy
between the constant path $H(0)(t) = \1$ and $H(1) = \Xi_\fast$.

Take $\alpha_1(s;t) = \tF_{\gamma_s}(t)$ and
$\alpha_{1/2}(s;t) = \alpha_0(s;t) \Xi_\fast(t)$.
Clearly,
$\alpha_1: \Ss^k \to \hat\Little\Ss^n(z) \subset \Omega\Spin_{n+1}(z)$,
as required. It suffices to construct homotopies between
$\alpha_0$ and $\alpha_{1/2}$ and between $\alpha_{1/2}$ and $\alpha_1$.
The homotopy between $\alpha_0$ and $\alpha_{1/2}$ is given by
\[ \alpha_\sigma(s;t) = \alpha_0(s;t) H(2\sigma)(t), \qquad
\sigma \in [0,1/2]. \]
Recall that $d(\alpha_{1/2}(s;t), \alpha_1(s;t)) < r_{n+1}/2$:
the homotopy between $\alpha_{1/2}$ and $\alpha_1$
is defined by joining these two points of $\Spin_{n+1}$
by the uniquely defined shortest geodesic
(parametrized by a constant multiple of arc length):
\[ \alpha_\sigma(s;t) = g_{\alpha_{1/2}(s;t),\alpha_1(s;t)}(2\sigma-1),
\qquad \sigma \in [1/2,1]. \]
\end{proof}

\begin{prop}
\label{prop:injective}
Assume $\Pi(z) = \pm J_+$:
then the inclusion $\hat\Little\Ss^n(z) \subset \Omega\Spin_{n+1}(z)$
is a weak homotopy equivalence.
In other words (given Proposition \ref{prop:surjective}),
if $\hat H_0: \BB^{k+1} \to \Omega\Spin_{n+1}(z)$
takes $\Ss^k \subset \BB^{k+1}$
to $\hat\Little\Ss^n(z) \subset \Omega\Spin_{n+1}(z)$
then there exist
$H_1: \BB^{k+1} \to \Little\Ss^n(z)$ and corresponding
$\hat H_1: \BB^{k+1} \to \hat\Little\Ss^n(z)$ with
$\hat H_0|_{\Ss^k} = \hat H_1|_{\Ss^k}$.
\end{prop}

Clearly if $\Pi(z) = \pm J_+$ then $s(z)$ must be $0$, $1$ or $-1$.
Theorem \ref{theo:dois} therefore follows
from Proposition \ref{prop:injective} and Fact \ref{fact:BH}.

\begin{proof}
Assume without loss of generality that $H_0$ is smooth.
Take  $\xi_\fast \in \Little\Ss^n(\1)$ as in Lemma \ref{lemma:fast}
so that, for any $s \in \BB^{k+1}$,
$\gamma(s) = \hat H_0(s) \xi_\fast \in \Little\Ss^n(z)$.
We may furthermore assume that the curves
$\hat H_0(s)(C_1 t) \xi_\fast(C_2 t + C_3)$
are locally convex for any $s \in \BB^{k+1}$,
for any $C_1, C_2 \in [1/10,10]$ and for any $C_3 \in \RR$.
Recall that $\xi_\fast(t) = \xi_1(Nt)$ for some large $N$:
take $N$ to be a multiple of $4$ so that
$\Xi_\fast(1/4) = \Xi_\fast(1/2) = \Xi_\fast(3/4) = \1$,
$\Xi_\fast(t) = \TR(\Xi_\fast(t)) = J_+ (\Xi_\fast(t))^{-1} J_+$
and $\Xi_\fast(1-t) = (\Xi_\fast(t))^{-1}$.
Recall that $\Lambda_{\xi_\fast}$ is constant:
let $B = \Lambda_{\xi_\fast}(t)$.
Set
\[ H_1(s)(t) = \gamma(2s)(t) = \hat H_0(s)(t) \xi_\fast(t),
\qquad |s| \le 1/2. \]
We now define $H_1$ in the two regions
$|s| \in [1/2,3/4]$ and $|s| \in [3/4,1]$.

For $s \in [3/4,1]$ we squeeze the function $\hat H_0(s/|s|)$
to a central interval $[1-|s|,|s|]$ and attach chunks of $\xi_\fast(2t)$
outside the central interval.
For 
\[ n = 2, \quad
Q = \begin{pmatrix} 0 & 0 & -1 \\ 0 & -1 & 0 \\ -1 & 0 & 0 \end{pmatrix}
= \chop(-J_+), \]
the construction is illustrated in
Figure \ref{fig:loopatendpoints}:
we can add chunks of a locally convex curve at the endpoints
and translate in the sphere (i.e., rotate in $\RR^3$)
the central portion of the curve.
Continue the process to add several closed circles
at both endpoints.

\begin{figure}[ht]
\begin{center}
\epsfig{height=50mm,file=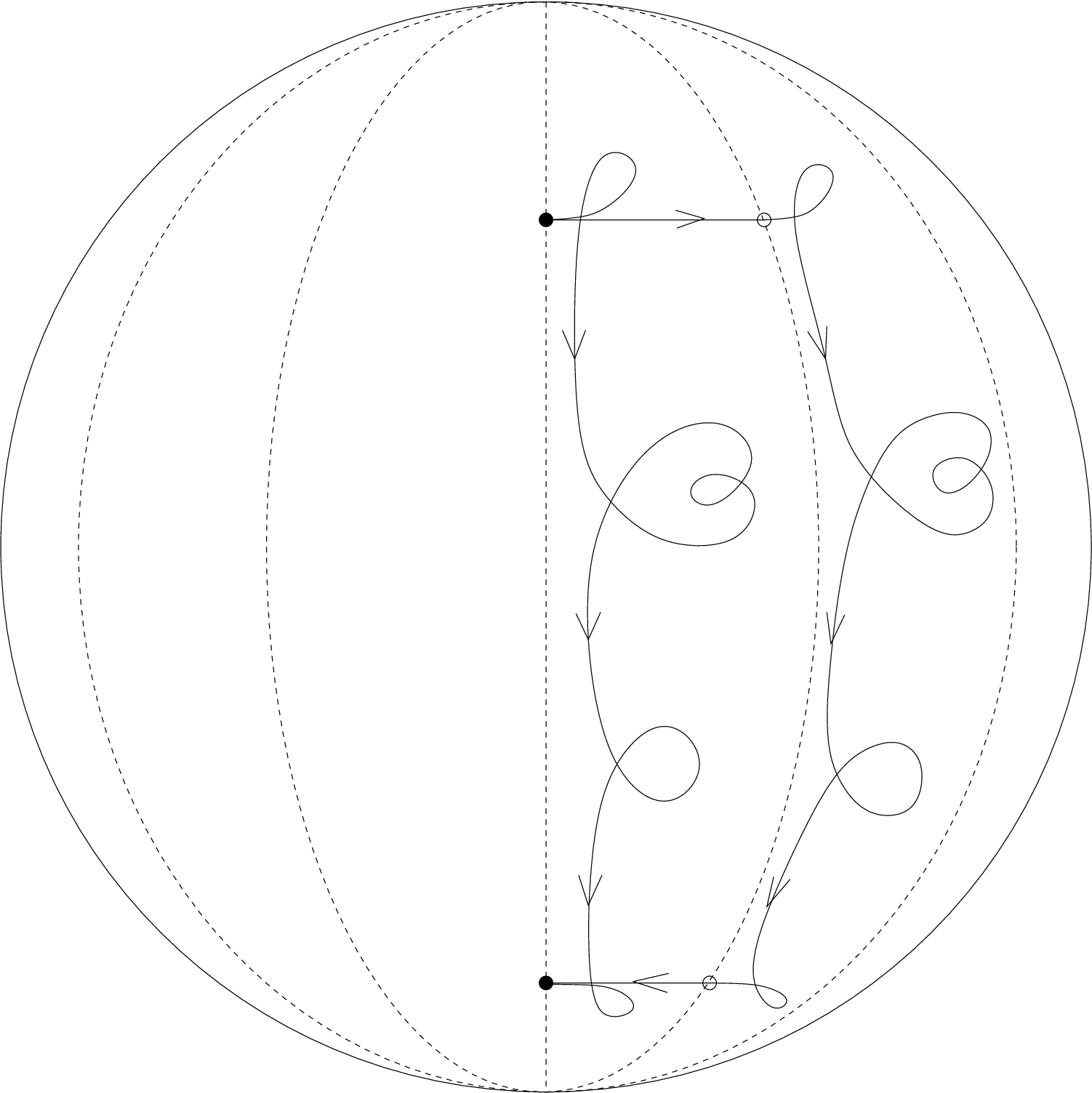}
\end{center}
\caption{Approximating a curve by a locally convex curve}
\label{fig:loopatendpoints}
\end{figure}

The general construction is perhaps best stated
in terms of $\Lambda$: for $|s| \in [3/4,1]$
\[ \Lambda_{H_1(s)}(t) = \begin{cases}
2B, & 0 \le t < 1-|s|, \\
\frac{1}{2|s| - 1}
\Lambda_{\hat H_0(s/|s|)}\left( \frac{t - 1 + |s|}{2|s| - 1} \right),
& 1-|s| \le t \le |s|, \\
2B, & |s| < t \le 1. \end{cases} \]
Recall that $\Lambda$ is only assumed to be of class $L^2$
and therefore the jump discontinuities are allowed. 
The curve $H_1(s)$ defined using the above $\Lambda$
is by construction locally convex:
we must verify that $\F_{H_1(s)}(1) = \hat H_1(s)(1) = z$.
We have $\hat H_1(s)(1-|s|) = \Xi_\fast(2(1-|s|))$;
for $1-|s| \le t \le |s|$ we therefore have
\[ \hat H_1(s)(t) =
\Xi_\fast(2(1-|s|)) {\hat H_0(s/|s|)}
\left( \frac{t - 1 + |s|}{2|s| - 1} \right) \]
and $\hat H_1(s)(|s|) = \Xi_\fast(2(1-|s|)) z$;
finally, at least in $SO_{n+1}$ we have
\begin{align*}
\hat H_1(s)(1) &= \Xi_\fast(2(1-|s|)) z  \Xi_\fast(2(1-|s|)) \\
&= \Xi_\fast(2(1-|s|)) (\pm J_+)  \Xi_\fast(2(1-|s|)) (\pm J_+) z \\
&= \Xi_\fast(2(1-|s|)) (\Xi_\fast(2(1-|s|)))^{-1} z = z
\end{align*}
(recall that $z = \pm J_+$ and that 
$J_+  \Xi_\fast(t) J_+ = (\Xi_\fast(t))^{-1}$);
by continuity we have $\hat H_1(s)(1) = z$ in $\Spin_{n+1}$
for all $s$ with $|s| \in [3/4,1]$.

The missing step is $|s| \in [1/2,3/4]$.
For $n = 2$, the circles which are concentrated at the endpoints
for $|s| = 3/4$ must spread along the curve as $s$ approaches $1/2$.
More algebraically, notice that
for both $|s| = 1/2$ and $|s| = 3/4$, we can write
$H_1(s)(t) = (A_{|s|}(s/|s|)(t))(\beta_{|s|}(t))$,
$A_{|s|}(s/|s|): [0,1] \to SO_{n+1}$, $\beta_{|s|}: [0,1] \to \Ss^n$.
From the constructions above we have
\begin{gather*}
A_{\frac12}(s/|s|)(t) = \hat H_0(s/|s|)(t), \qquad
\beta_{\frac12}(t) = \xi_\fast(t), \\
A_{\frac34}(s/|s|)(t) = \hat H_0(s/|s|)(g_{\frac34}(t)), \qquad
\beta_{\frac34}(t) = \xi_\fast(h_{\frac34}(t)) 
\end{gather*}
where
\[
g_{\frac{3}{4}}(t) = \begin{cases} 0,& 0 \le t \le \frac{1}{4}, \\
2t-\frac{1}{2}, & \frac{1}{4} \le t \le \frac{3}{4}, \\
1, & \frac{3}{4} \le t \le 1, \end{cases} \\
\quad
h_{\frac{3}{4}}(t) = \begin{cases} 2t,& 0 \le t \le \frac{1}{4}, \\
\frac{1}{2}, & \frac{1}{4} \le t \le \frac{3}{4}, \\
2t - 1, & \frac{3}{4} \le t \le 1. \end{cases} \]
We complete the definition of $H_1$ with
\begin{gather*}
H_1(s)(t) = (A_{|s|}(s/|s|)(t))(\beta_{|s|}(t)), \qquad 1/2 \le |s| \le 3/4 \\
A_{\sigma}(s/|s|)(t) = \hat H_0(s/|s|)(g_{\sigma}(t)), \qquad
\beta_{\sigma}(t) = \xi_\fast(h_{\sigma}(t)),
\end{gather*}
$g_\sigma$ and $h_\sigma$ as plotted in Figure \ref{fig:ghsigma}
(notice that $g_{\frac12}(t) = h_{\frac12}(t) = t$).

\begin{figure}[ht]
\psfrag{g}{$g$}
\psfrag{h}{$h$}
\psfrag{1}{$1$}
\psfrag{1/2}{$\frac12$}
\psfrag{s}{$\sigma$}
\psfrag{1-s}{$1-\sigma$}
\psfrag{s-1/2}{$\sigma-\frac12$}
\psfrag{3/2-s}{$\frac32-\sigma$}
\begin{center}
\epsfig{height=40mm,file=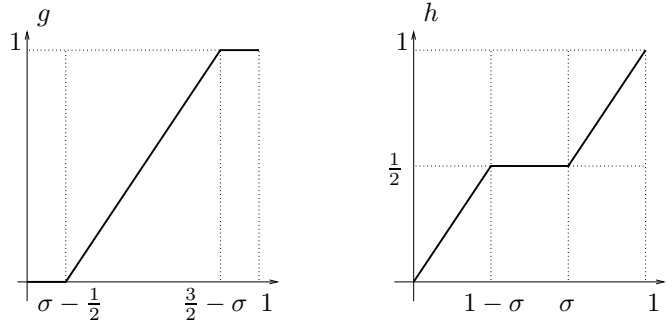}
\end{center}
\caption{The functions $g_\sigma$ and $h_\sigma$}
\label{fig:ghsigma}
\end{figure}

We are left with proving that $H_1(s)$ is locally convex.
For $t \in [0,|s|-\frac12] \cup [\frac32-|s|,1]$,
$H_1(s)$ is a reparametrization of $\xi_\fast$ and therefore locally convex.
For $t \in [|s|-\frac12,1-|s|] \cup [|s|,\frac32-|s|]$,
locally convexity follows from Lemma \ref{lemma:fast} or,
perhaps more precisely, from the choice of $\xi_\fast$
as described at the beginning of the proof.
Finally, for $t \in [1-|s|,|s|]$,
$H_1(s)$ is a reparametrization of $H_0(s/|s|)$
and therefore again locally convex.
This completes the construction of $H_1$ and the proof.
\end{proof}

\bigbreak

\section{Final remarks and open problems} 

\subsection{Is Theorem \ref{theo:main} strong?}
For $n = 2$, Theorems \ref{theo:main} and \ref{theo:dois} imply
that any space $\Little\Ss^2(z)$ is homeomorphic to one of three spaces 
$\Little\Ss^2(\1)$, $\Little\Ss^2(-\1)$
or $\Omega\Spin(3) = \Omega\Ss^3$.
From \cite{Little}, we know that $\Little\Ss^2(\1)$ and $\Little\Ss^2(-\1)$
have $1$ and $2$ connected components, respectively,
and $\Omega\Ss^3$ is clearly connected.
From \cite{Saldanha1} and \cite{Saldanha2}, we know that
$\dim H^2(\Little\Ss^2(\1);\RR) = 2$,
$\dim H^2(\Little\Ss^2(-\1);\RR) = 1$
and $\dim H^4(\Little\Ss^2(-\1);\RR) \ge 2$.
Thus, these three spaces are not pairwise homeomorphic;
also, the non-contractible connected component of
$\Little\Ss^2(-\1)$ is not homeomorphic to
either $\Omega\Ss^3$ or $\Little\Ss^2(-\1)$.

Unfortunately, similar information is unavailable for $n>2$.
We formulate the following question. 

\begin{prob}
\label{prob:mainppisstrong}
Are the $\lceil\frac{n}{2}\rceil+1$ subspaces
$\Little\Ss^n(M^{n+1}_s)$
(and similar space of curves in $\Spin_{n+1}$)
appearing in Theorem~\ref{theo:main}
pairwise non-homeomorphic for $n > 2$? 
\end{prob} 

Our best guess is that the answer is positive.

\subsection{Bounded curvature}
A first natural generalization of the space of locally curves on
$\Ss^2$ is the space of curves whose curvature $\kappa$ at each point is
bounded by two constants $m < \kappa <M$. 

\begin{prob}
Is it true that there are only
finitely many topologically distinct spaces of curves whose curvature is
bounded as above among the spaces of such curves with the fixed initial and
variable finite frames?  
\end{prob}

\subsection{Other Lie groups}
The space $\hat\Little\Ss^n$ is a special instance
of a more general construction on an arbitrary compact Lie group.
Given a compact Lie group $G$, consider a non-holonomic subspace 
of its Lie algebra (i.e., this subspace generates the whole algebra). 
Consider some polytopal convex cone in this subspace.
Take the left-invariant distribution of cones on $G$ obtained
by its left translation  in the algebra.
Finally, consider spaces of curves on $G$ tangent to the
obtained cone distribution which start at the unit element and end at some
fixed point of $G$. 

This generalization includes the scenario described in the previous
subsection as a special case ($G$ is $SO_{3}$ and the subspace consists
of skew tridiagonal matrices, just as for our problem;
the only difference is the cone).

\begin{prob}
Is it true that there are only
finitely many topologically distinct spaces of such curves with the fixed
initial and variable finite point?  
\end{prob}

This is likely to be too optimistic an attempt of  generalization,
but perhaps the finiteness condition holds true
with some interesting additional hypothesis.
For instance, our cone is the interior of the convex hull of a small
set of rather special vectors: maybe some such condition is needed.

\subsection{The homotopy type of spaces of closed locally convex curves}
Finally, the most interesting problem in this context
is to describe the homotopy type of the space of closed locally convex curves.
The aim of \cite{Saldanha3} is to address this problem for $n = 2$;
see partial results in \cite{Saldanha1}, \cite{Saldanha2}.


\bigbreak

\bigbreak

{

\parindent=0pt
\parskip=0pt
\obeylines

Nicolau C. Saldanha, Departamento de Matem\'atica, PUC-Rio
R. Marqu\^es de S. Vicente 225, Rio de Janeiro, RJ 22453-900, Brazil
saldanha@puc-rio.br; http://www.mat.puc-rio.br/$\sim$nicolau

\smallskip

Boris Shapiro, Stockholm University, S-10691, Stockholm, Sweden
shapiro@math.su.se; http://www.math.su.se/$\sim$shapiro

}

\end{document}